\documentclass[12pt,a4paper]{article}
\usepackage{amsfonts}
\usepackage{amssymb}
 \newcommand{\R}{\mathbb{R}}
\newcommand{\C}{\mathbb{C}}

\newcommand{\N}{\mathbb{N}}
\newtheorem{theo}{Theorem}[section]
\newtheorem{lemm}{Lemma}[section]
\newtheorem{prop}{Proposition}[section]
\newtheorem{coro}{Corollary}[section]
\newtheorem{rema}{Remark}[section]
\newcommand{\dint}{\displaystyle\int}
\newcommand{\dfrac}{\displaystyle\frac}

\newcommand{\norm}[1]{\left\Vert#1\right\Vert}
\newcommand{\abs}[1]{\left\vert#1\right\vert}
\title{\bf Logarithmic decay of the energy for an hyperbolic-parabolic coupled system}
\author{{\bf Ines Kamoun Fathallah}\thanks{Laboratoire LMV, Universit\'e de Versailles
Saint-Quentin-En-Yvelines, 45 Avenue des Etats-Unis Batiment Fermat
78035 Versailles (France). E-mail address:
ines.fathallah@math.uvsq.fr; Tel:+33139253629; Fax:+33139254645.}}
\date{}
\begin{document}
\maketitle
\begin{abstract}
This paper is devoted to the study of a coupled system consisting in
a wave and heat equations coupled through transmission condition
along a steady interface. This system is a linearized model for
fluid-structure interaction introduced by Rauch, Zhang and Zuazua
for a simple transmission condition and by Zhang and Zuazua for a
natural transmission condition.\par Using an abstract Theorem of
Burq and a new Carleman estimate shown near the interface, we
complete the results obtained by Zhang and Zuazua and by Duyckaerts.
We show, without any geometric restriction, a logarithmic decay
result.
\end{abstract}

\section{Introduction and results}
\hspace{5mm}In this work, we are interested with a linearized model
for fluid-structure interaction introduced by Zhang and Zuazua in
\cite{ZZ} and Duyckaerts in \cite{TD}. This model consists of a wave
and heat equations coupled through an interface with suitable
transmission conditions. Our purpose is to analyze the stability of
this system and so to determine the decay rate of energy of solution
as $t\rightarrow\infty$.
\par Let $\Omega\subset \R^{n}$ be a bounded
domain with a smooth boundary $\Gamma=\partial \Omega$. Let
$\Omega_{1}$ and $\Omega_{2}$ be two bounded open sets with smooth
boundary such that $\Omega_{1}\subset\Omega$ and
$\Omega_{2}=\Omega\backslash \overline {\Omega}_{1}$. We denote by
$\gamma= \partial \Omega_{1}\cap
\partial\Omega_{2}$ the interface, $ \gamma\subset\subset \Omega$, $\Gamma_{j}=\partial
\Omega_{j}\backslash\gamma$, $j=1,2$, $\partial_{n}$ and
$\partial_{n'}$ the unit outward normal vectors of $\Omega_{1}$ and
$\Omega_{2}$ respectively ($\partial_{n'}=-\partial_{n}$ on
$\gamma$).
\begin{equation}
\left \{\begin{array}{ll}\label{s1}
\partial_{t}u-\triangle u=0&\mbox{in}\,(0,\infty)\times \Omega_{1},\\
\partial_{t}^{2}v-\triangle v=0&\mbox{in}\,(0,\infty)\times \Omega_{2},\\
u=0&\mbox{on}\,(0,\infty)\times\Gamma_{1},\\
v=0&\mbox{on}\,(0,\infty)\times\Gamma_{2},\\
u=\partial_{t}v,\quad
\partial_{n}u=-\partial_{n'}v&\mbox{on}\,(0,\infty)\times \gamma,\\
u|_{t=0}=u_{0}\in L^{2}(\Omega_{1})&\mbox{in}\,\Omega_{1},\\
v|_{t=0}=v_{0}\in H^{1}(\Omega_{2}), \quad
\partial_{t}v|_{t=0}=v_{1}\in L^{2}(\Omega_{2})&\mbox{in}\,\Omega_{2}.
\end{array} \right.
\end{equation}
In this system, $u$ may be viewed as the velocity of fluid; while
$v$ and $\partial_{t}v$ represent respectively the displacement and
velocity of the structure. That's why the transmission condition
$u=\partial_{t}v$ is considered as the natural condition. For
 the modelisation subject, we refer to \cite{RZZ} and \cite{ZZ}.
\par System (\ref{s1}) is introduced by Zhang and Zuazua \cite{ZZ}.
The same system was considered by Rauch, Zhang and Zuazua
 in \cite{RZZ} but for simplified
transmission condition $u=v$ on the interface instead of
$u=\partial_{t}v$. They prove, under a suitable Geometric Control
Condition (GCC) (see \cite{BLR}), a polynomial decay result. Zhang
and Zuazua in \cite{ZZ} prove, without GCC, a logarithmic decay
result. Duyckaerts in \cite{TD} improves these results.\par For
system (\ref{s1}), Zhang and Zuazua in \cite{ZZ}, show the lack of
uniform decay and they prove, under GCC, a polynomial decay result.
Without geometric conditions, they analyze the difficulty to prove
the logarithmic decay result. This difficulty is mainly due to the
lack of gain regularity of wave component $v$ near the interface
$\gamma$ (see \cite{ZZ}, Remark 19) which means that the embedding
of the domain $D(\mathcal{A})$ of dissipative operator in the energy
space is not compact (see \cite{ZZ}, Theorem 1). In \cite{TD},
Duyckaerts improves the polynomial decay result under GCC and
confirms the same obstacle to show the logarithmic decay for
solution of (\ref{s1}) without GCC. In this paper we are interested
with this problem.
\par There is an extensive literature on the stabilization of PDEs
and on the Logarithmic decay of the energy (\cite{MB2}, \cite{MB}
\cite{MB1},
 \cite{L}, \cite{LR2}, \cite{LR1} and the references cited therein)
 and this paper use a part of the idea developed in \cite{MB}.
 \par Here we recall the mathematical frame work for this problem (see
 \cite{ZZ}).
\par Define the energy space $H$ and the operator
$\mathcal{A}$ on $H$, of domain $D(\mathcal{A})$ by $$ H=\left\{
U_{0}=(u_{0},v_{0},v_{1})\in L^{2}(\Omega_{1})\times
H^{1}_{\Gamma_{2}}(\Omega_{2})\times L^{2}(\Omega_{2})\right\}$$
when $H^{1}_{\Gamma_{2}}(\Omega_{2}) $ is defined as the space
$$ H^{1}_{\Gamma_{2}}(\Omega_{2})= \left\{ v_{0}\in H^{1}(\Omega_{2}), v_{0}|_{\Gamma_{2}}=0\right\},$$
$$ \mathcal{A}U_{0}=( \triangle u_{0},v_{1}, \triangle v_{0})$$
\begin{eqnarray*}D(\mathcal{A})=\{ U_{0}\in H,\, u_{0}\in
H^{1}(\Omega_{1}),\, \triangle u_{0}\in L^{2}(\Omega_{1}),\,
~~~~~~~~~~~~~~~~~~~~~~~~~~~~~~~~~~~~~~\\\\
v_{1}\in H^{1}_{\Gamma_{2}}(\Omega_{2}),\, \triangle v_{0}\in
L^{2}(\Omega_{2}),\,
u_{0}|_{\gamma}=v_{1}|_{\gamma},\,\partial_{n}u_{0}|_{\gamma}=-\partial
_{n}v_{0}|_{\gamma}\}.\end{eqnarray*} System (\ref{s1}) may thus be
rewritten in the abstract form
$$ \partial_{t}U= \mathcal{A}U,\quad\quad U(t)=(u(t),v(t),\partial_{t}v(t)).$$
\hspace{5mm} For any solution $(u,v,\partial_{t}v)$ of system
(\ref{s1}), we have a natural energy
$$ E(t)=E(u,v,\partial_{t}v)(t)=\frac{1}{2}\left( \int_{\Omega_{1}}\abs{u(t)}^{2}dx+
\int_{\Omega_{2}}\abs{\partial_{t}v(t)}^{2}dx+
\int_{\Omega_{2}}\abs{\nabla v(t)}^{2}dx \right).$$ By means of the
classical energy method, we have
$$\frac{d}{dt} E(t)=-\int_{\Omega_{1}}\abs{\nabla u}^{2}dx.$$
Therefore the energy of (\ref{s1}) is decreasing with respect to
$t$, the dissipation coming from the heat component $u$. Our main
goal is to prove a logarithmic decay without any geometric
restrictions.\par As Duyckaerts \cite{TD} did for the simplified
model, the idea is, first, to use a known result of Burq (see
\cite{NB}) which links, for dissipative operators, logarithmic decay
to resolvent estimates with exponential loss; secondly to prove,
following the work of Bellassoued in \cite{MB}, a new Carleman
inequality near the interface $\gamma$.
\par The main results are given by Theorem \ref{TR} for resolvent and
Theorem \ref{DE} for decay.
\begin{theo}\label{TR}
There exists $C>0$, such that for every $\mu \in \R$ with
$\abs{\mu}$ large,\\ we have
\begin{equation}\label{R}
\norm{(\mathcal{A}-i\mu)^{-1}}_{\mathcal{L}(H)}\leq Ce^{C\abs{\mu}}.
\end{equation}
\end{theo}
\begin{theo}\label{DE}
There exists $C>0$, such that the energy of a smooth solution of
(\ref{s1}) decays at logarithmic speed
\begin{equation}\label{d}
\sqrt{E(t)}\leq \frac{C}{\log(t+2)}\norm{U}_{D(\mathcal{A})}.
\end{equation}
\end{theo}
\hspace{5mm} Burq in (\cite{NB}, Theorem 3) and Duyckaerts in
(\cite{TD}, Section 7) show that to prove Theorem \ref{DE} it
suffices to show Theorem \ref{TR}.
\par The strategy of the proof of
Theorem \ref{TR} is the following. A new Carleman estimate shown
near the interface $\gamma$ implies an interpolation inequality
given by Theorem \ref{II}. Theorem \ref{II} implies Theorem \ref{TI}
which gives an estimate of the wave component by the heat one and
which is the key point of the proof of Theorem \ref{TR}.\par The
rest of this paper is organized as follows. In section 2, we show,
from Theorem \ref{TI}, Theorem \ref{TR} and we explain how Theorem
\ref{II} implies Theorem \ref{TI}. In section 3, we begin by stating
the new Carleman estimate and explain how this estimate implies
Theorem \ref{II}. We give then the proof of this Carleman estimate.
Section 4 is devoted to the proof of important estimates stated in
Theorem \ref{t2} in the proof of Carleman estimate. Appendices A and
B are devoted to prove some technical results that will be used
along the paper.\\\\
{\bf{Acknowledgment}}\\\\
Sincere thanks to professor Luc Robbiano for inspiring question, his
greatly contribution to this work and for his careful reading of the
manuscript. I want to thank also professor Mourad Bellassoued for
his proposition to work in this domain.
\section{Proof of Theorem \ref{TR}}
\hspace{5mm} This section is devoted to the proof of Theorem
\ref{TR}. We start by stating Theorem \ref{TI}. Then we will explain
how this Theorem implies Theorem \ref{TR}. Finally, we give the
proof of Theorem \ref{TI}.\par
 Let $\mu$ be a real number such that
$\abs{\mu}$ is large, and assume
\begin{equation}\label{s2}
F=(\mathcal{A}-i\mu)U, \quad U=(u_{0}, v_{0}, v_{1})\in
D(\mathcal{A}),\quad\quad F=(f_{0},g_{0}, g_{1})\in H
\end{equation}
The equation (\ref{s2}) yields \begin{equation} \left
\{\begin{array}{rclc}\label{s12}
(\triangle-i\mu) u_{0}&=&f_{0}& \mbox{in}\, \Omega_{1},\\
(\triangle+\mu^{2}) v_{0}&=&g_{1}+i\mu g_{0}& \mbox{in}\, \Omega_{2},\\
v_{1}&=&g_{0}+i\mu v_{0}&\mbox{in}\, \Omega_{2},
\end{array} \right.\end{equation}
with the following boundary conditions \begin{equation}\label{bis}
\left \{\begin{array}{lcl}
u_{0}|_{\Gamma_{1}}=0,\,\, v_{0}|_{\Gamma_{2}}=0 &&\\
op(b_{1})u=u_{0}-i\mu v_{0}&=&g_{0}|_{\gamma},\\
op(b_{2})u=\partial_{n}u_{0}-\partial_{n}v_{0}&=&0|_{\gamma}.
\end{array} \right.\end{equation}
To proof Theorem \ref{TR}, we begin by stating this result
\begin{theo}\label{TI}
Let $U=(u_{0}, v_{0}, v_{1})\in D(\mathcal{A})$ satisfying equation
(\ref{s12}) and (\ref{bis}). Then there exists constants $C>0$,
$c_{1}>0$ and $\mu_{0}>0$ such that for any $\mu\geq\mu_{0}$ we have
the following estimate
\begin{equation}\label{s11}
\norm{v_{0}}^{2}_{H^{1}(\Omega_{2})}\leq
Ce^{c_{1}\mu}\left(\norm{f_{0}}^{2}_{L^{2}(\Omega_{1})}+
 \norm{ g_{1}+i\mu
g_{0}}^{2}_{L^{2}(\Omega_{2})}+\norm{g_{0}}^{2}_{H^{1}(\Omega_{2})}+
\norm{u_{0}}^{2}_{H^{1}(\Omega_{1})} \right).
\end{equation}
\end{theo}
Moreover, from the first equation of system (\ref{s12}), we have
$$\int_{\Omega_{1}}(-\triangle +i\mu)u_{0}\overline{u_{0}}dx=\norm{\nabla u_{0}}_{L^{2}(\Omega_{1})}^{2}+i\mu
\norm{u_{0}}^{2}_{L^{2}(\Omega_{1})}-\int_{\gamma}\partial_{n}u_{0}\overline{u_{0}}d\sigma
.$$Since $u_{0}|_{\gamma}=g_{0}+i\mu v_{0}$ and
$\partial_{n}u_{0}=-\partial_{n'}v_{0}$, then
\begin{eqnarray}\label{s13}
\int_{\Omega_{1}}(-\triangle
+i\mu)u_{0}\overline{u}_{0}dx=\norm{\nabla
u_{0}}_{L^{2}(\Omega_{1})}^{2}+i\mu
\norm{u_{0}}^{2}_{L^{2}(\Omega_{1})}-i\mu\int_{\gamma}\partial_{n'}v_{0}\overline{v}_{0}d\sigma
+\int_{\gamma}\partial_{n'}v_{0}\overline{g}_{0}d\sigma.
\end{eqnarray}
From the second equation of system (\ref{s12}) and multiplying by
$(-i\mu)$, we obtain
\begin{equation}\label{s14}
i\mu\int_{\Omega_{2}}(\triangle+\mu^{2})v_{0}\overline{v}_{0}dx=-i\mu\norm{\nabla
v_{0}}_{L^{2}(\Omega_{2})}^{2}+i\mu^{3}
\norm{v_{0}}^{2}_{L^{2}(\Omega_{2})}+i\mu\int_{\gamma}\partial_{n'}v_{0}\overline{v}_{0}d\sigma.
\end{equation}
Adding (\ref{s13}) and (\ref{s14}), we obtain
\begin{eqnarray*}
\int_{\Omega_{1}}(-\triangle
+i\mu)u_{0}\overline{u}_{0}dx+i\mu\int_{\Omega_{2}}(\triangle+\mu^{2})v_{0}\overline{v}_{0}dx=
~~~~~~~~~~~~~~~~~~~~\\
i\mu \norm{u_{0}}^{2}_{L^{2}(\Omega_{1})}+\norm{\nabla
u_{0}}_{L^{2}(\Omega_{1})}^{2}-i\mu\norm{\nabla
v_{0}}_{L^{2}(\Omega_{2})}^{2}+i\mu^{3}
\norm{v_{0}}^{2}_{L^{2}(\Omega_{2})}+\int_{\gamma}\partial_{n'}v_{0}\overline{g}_{0}d\sigma.
\end{eqnarray*}
Taking the real part of this expression, we get
\begin{equation}\label{s15}
\norm{\nabla u_{0}}_{L^{2}(\Omega_{1})}^{2}\leq \norm{ (\triangle
-i\mu)u_{0}}_{L^{2}(\Omega_{1})}\norm{u_{0}}_{L^{2}(\Omega_{1})}+
\norm{(\triangle+\mu^{2})v_{0}}_{L^{2}(\Omega_{2})}\norm{v_{0}}_{L^{2}(\Omega_{2})}+
\abs{\int_{\gamma}\partial_{n'}v_{0} \overline{g}_{0}d\sigma}.
\end{equation}
Recalling that $ \triangle v_{0}=g_{0}+i\mu g_{0}-\mu^{2}v_{0}$ and
using the trace lemma (Lemma 3.4 in \cite{TD}), we obtain
$$ \norm{\partial_{n}v_{0}}_{H^{-\frac{1}{2}}(\gamma)}\leq C\left( \mu^{2}\norm{v_{0}}_{H^{1}(\Omega_{2})}
+ \norm{ g_{1}+i\mu g_{0}}_{L^{2}(\Omega_{2})}\right).$$ Combining
with (\ref{s15}), we obtain
\begin{eqnarray*}
 \norm{\nabla u_{0}}_{L^{2}(\Omega_{1})}^{2}\leq \norm{f_{0}}_{L^{2}(\Omega_{1})}\norm{u_{0}}_{L^{2}(\Omega_{1})}+
\norm{ g_{1}+i\mu g_{0}}_{L^{2}(\Omega_{2})}\norm{v_{0}}_{L^{2}(\Omega_{2})}  ~~~~~~~~~~~~~~
\nonumber\\+ \left(
\mu^{2}\norm{v_{0}}_{H^{1}(\Omega_{2})}+ \norm{ g_{1}+i\mu
g_{0}}_{L^{2}(\Omega_{2})}\right)\norm{g_{0}}_{H^{\frac{1}{2}}(\gamma)}.
\end{eqnarray*}
Then \begin{eqnarray*}\norm{\nabla u_{0}}_{L^{2}(\Omega_{1})}^{2}\leq
\frac{C}{\epsilon} \norm{f_{0}}_{L^{2}(\Omega_{1})}^{2}+
\epsilon\norm{u_{0}}_{L^{2}(\Omega_{1})}^{2}+
 \frac{C}{\epsilon} \norm{ g_{1}+i\mu g_{0}}_{L^{2}(\Omega_{2})}^{2}+
 \epsilon \norm{v_{0}}^{2}_{L^{2}(\Omega_{2})}~~~~~~~~
 \nonumber\\+
 \left(\mu^{2}
\norm{v_{0}}_{H^{1}(\Omega_{2})}+ \norm{ g_{1}+i\mu
g_{0}}_{L^{2}(\Omega_{2})}\right)\norm{g_{0}}_{H^{\frac{1}{2}}(\gamma)}.
\end{eqnarray*}
Now we need to use this result shown in Appendix A.
\begin{lemm}\label{le1}
Let $\mathcal{O}$ be a bounded open set of $\R^{n}$. Then there
exists $C>0$ such that for $u$ and $f$ satisfying $
(\triangle-i\mu)u=f$ in $ \mathcal{O}$, $\mu\geq 1$, we have the
following estimate
\begin{equation}\label{le11}\norm{u}_{H^{1}(\mathcal{O})}\leq
C\left(\norm{\nabla u}_{L^{2}(\mathcal{O})}+
\norm{f}_{L^{2}(\mathcal{O})}\right).\end{equation}
\end{lemm}
Using this Lemma, we obtain, for $\epsilon$ small enough
\begin{eqnarray*} \norm{u_{0}}_{H^{1}(\Omega_{1})}^{2}\leq C \norm{f_{0}}_{L^{2}(\Omega_{1})}^{2}+
C_{\epsilon} \norm{ g_{1}+i\mu g_{0}}_{L^{2}(\Omega_{2})}^{2}+
\epsilon \norm{v_{0}}^{2}_{L^{2}(\Omega_{2})}
~~~~~~~~~~\nonumber\\+ \left(\mu^{2} \norm{v_{0}}_{H^{1}(\Omega_{2})}+ \norm{ g_{1}+i\mu
g_{0}}_{L^{2}(\Omega_{2})}\right)\norm{g_{0}}_{H^{\frac{1}{2}}(\gamma)}.
\end{eqnarray*}
Then there exists $c_{3}>>c_{1}$ such that
\begin{equation}\label{s16}
\norm{u_{0}}^{2}_{H^{1}(\Omega_{1})}\leq
C\left(\norm{f_{0}}_{L^{2}(\Omega_{1})}^{2}+
\epsilon e^{-c_{3}\mu}\norm{v_{0}}_{H^{1}(\Omega_{2})}^{2}+C_{\epsilon}
e^{-c_{3}\mu}\norm{ g_{1}+i\mu
g_{0}}_{L^{2}(\Omega_{2})}^{2}+e^{c_{3}\mu}\norm{g_{0}}_{H^{1}(\Omega_{2})}^{2}
\right).
\end{equation}
Inserting in (\ref{s11}), we obtain, for $\epsilon$ small enough
\begin{equation}\label{s17}
\norm{v_{0}}^{2}_{H^{1}(\Omega_{2})}\leq
Ce^{c\mu}\left(\norm{f_{0}}_{L^{2}(\Omega_{1})}^{2}+\norm{g_{0}}_{H^{1}(\Omega_{2})}^{2}+
\norm{ g_{1}+i\mu g_{0}}_{L^{2}(\Omega_{2})}^{2}\right).
\end{equation}
Combining (\ref{s16}) and (\ref{s17}), we obtain
\begin{equation}\label{s18}
\norm{u_{0}}^{2}_{H^{1}(\Omega_{1})}\leq
Ce^{c\mu}\left(\norm{f_{0}}_{L^{2}(\Omega_{1})}^{2}+\norm{g_{0}}_{H^{1}(\Omega_{2})}^{2}+
\norm{ g_{1}+i\mu g_{0}}_{L^{2}(\Omega_{2})}^{2}\right).
\end{equation}
Recalling that $ v_{1}=g_{0}+i\mu v_{0}$ and using (\ref{s17}), we
obtain
\begin{equation}\label{s19}
\norm{v_{1}}^{2}_{H^{1}(\Omega_{1})}\leq
Ce^{c\mu}\left(\norm{f_{0}}_{L^{2}(\Omega_{1})}^{2}+\norm{g_{0}}_{H^{1}(\Omega_{2})}^{2}+
\norm{ g_{1}+i\mu g_{0}}_{L^{2}(\Omega_{2})}^{2}\right).
\end{equation}
Combining (\ref{s17}), (\ref{s18}) and (\ref{s19}), we obtain
Theorem
\ref{TR}.\\$~~~~~~~~~~~~~~~~~~~~~~~~~~~~~~~~~~~~~~~~~~~~~~~~~~~~~~~~~~~~~~~~~~~~~~~~~~~~~~~~~~~~~~
 ~~~~~~~~~~~~~~~~~~~~~~~~~~~~~~~~~~~~~~~~~~~~~~~~\square$
{\bf Proof of Theorem \ref{TI}}\\\\
Estimate (\ref{s11}) is consequence of two important results. The
first is a known result shown by Lebeau and Robbiano in \cite{LR}
and the second one is given by Theorem \ref{II} and proved in
section 3.\\ Let $0<\epsilon_{1}<\epsilon_{2}$ and
$V_{\epsilon_{j}}$, $j=1,2$, such that
$V_{\epsilon_{j}}=\{x\in\Omega_{2},\,\,
d(x,\gamma)<\epsilon_{j}\}$.\\
 Recalling that $(\triangle+\mu^{2})v_{0}=g_{1}+i\mu g_{0}$, then for all $D>0$, there
exists $C>0$ and $ \nu\in ]0,1[$ such that we have the following
estimate (see \cite{LR})
\begin{equation}\label{a1}
\norm{v_{0}}_{H^{1}(\Omega_{2}\backslash V_{\epsilon_{1}})}\leq
Ce^{D\mu}\norm{v_{0}}_{H^{1}(\Omega_{2})}^{1-\nu}\left(\norm{ g_{1}+i\mu
g_{0}}_{L^{2}(\Omega_{2})}+
\norm{v_{0}}_{H^{1}(V_{\epsilon_{2}})}\right)^{\nu}
\end{equation}
Moreover we have the following result shown in section 3.
\begin{theo}\label{II}
There exists $C>0$, $c_{1}>0$, $c_{2}>0$, $\epsilon_{2}>0$ and
$\mu_{0}>0$ such that for any $\mu\geq\mu_{0}$, we have the
following estimate
\begin{eqnarray}\label{II1}
\norm{v_{0}}^{2}_{H^{1}(V_{\epsilon_{2}})}&\leq& C
e^{c_{1}\mu}\left[\norm{f_{0}}^{2}_{L^{2}(\Omega_{1})}+
 \norm{ g_{1}+i\mu
g_{0}}^{2}_{L^{2}(\Omega_{2})}+\norm{g_{0}}^{2}_{H^{1}(\Omega_{2})}+
\norm{u_{0}}^{2}_{H^{1}(\Omega_{1})}\right]\nonumber\\\nonumber\\&+&C
e^{-c_{2}\mu}\norm{v_{0}}_{H^{1}(\Omega_{2})}^{2}.
\end{eqnarray}
\end{theo}
 Combining (\ref{a1})
and (\ref{II1}) we obtain
\begin{eqnarray}\label{a3}\norm{v_{0}}_{H^{1}(\Omega_{2}\backslash
V_{\epsilon_{2}})}^{2}\leq C
\epsilon e^{D\mu}\norm{v_{0}}_{H^{1}(\Omega_{2})}^{2}+
\frac{C}{\epsilon}\norm{ g_{1}+i\mu g_{0}}_{L^{2}(\Omega_{2})}^{2}+
 \frac{C}{\epsilon}e^{-c_{2}\mu}\norm{v_{0}}_{H^{1}(\Omega_{2})}^{2}~~~~\nonumber\\
+\frac{C}{\epsilon}e^{c_{1}\mu}\left[\norm{f_{0}}^{2}_{L^{2}(\Omega_{1})}
+ \norm{ g_{1}+i\mu
g_{0}}^{2}_{L^{2}(\Omega_{2})}+\norm{g_{0}}^{2}_{H^{1}(\Omega_{2})}+
\norm{u_{0}}^{2}_{H^{1}(\Omega_{1})}\right].
\end{eqnarray}
Adding (\ref{II1}) and (\ref{a3}), we obtain
\begin{eqnarray*}
\norm{v_{0}}^{2}_{H^{1}(\Omega_{2})}\leq
C\epsilon e^{D\mu} \norm{v_{0}}^{2}_{H^{1}(\Omega_{2})}+C_{\epsilon}\norm{
g_{1}+i\mu
g_{0}}^{2}_{L^{2}(\Omega_{2})}+C_{\epsilon}e^{-c_{2}\mu}\norm{v_{0}}_{H^{1}(\Omega_{2})}^{2}
~~~~~~~~~~~~~~~~~\\+C_{\epsilon}e^{c_{1}\mu}\left[\norm{f_{0}}^{2}_{L^{2}(\Omega_{1})}+
 \norm{ g_{1}+i\mu
g_{0}}^{2}_{L^{2}(\Omega_{2})}+\norm{g_{0}}^{2}_{H^{1}(\Omega_{2})}+
\norm{u_{0}}^{2}_{H^{1}(\Omega_{1})}\right].
\end{eqnarray*}
We fixe $ \epsilon$ small enough and $D<c_{2}$, then there exists $\mu_{0}>0$ such
that for any $\mu\geq \mu_{0}$, we obtain (\ref{s11}).\\
$~~~~~~~~~~~~~~~~~~~~~~~~~~~~~~~~~~~~~~~~~~~~~~~~~~~~~~~~~~~~~~~~~~~~~~~~~~~~~~~~~~~~~~~~~~~~~~~~
~~~~~~~~~~~~~~~~~~~~~~~~~~~~~~~~~~~~~~~~~~~~~~~~~~~~~~~~~~\square$
\section{Carleman estimate and Consequence}
\hspace{5mm} In this part, we show the new Carleman estimate and we
prove Theorem \ref{II} which is consequence of this estimate.
\subsection{State of Carleman estimate}
\hspace{5mm} In this subsection we state the Carleman estimate which
is the starting point of the proof of the main result. Let
$u=(u_{0}, v_{0})$ satisfies the equation
\begin{equation}
\left \{\begin{array}{lc}\label{e1}
-(\triangle+\mu)u_{0}=f_{1}&\mbox{in}\, \Omega_{1},\\
-(\triangle+\mu^{2})v_{0}=f_{2}&\mbox{in}\, \Omega_{2},\\
op(B_{1})u=u_{0}-i\mu v_{0}=e_{1}&\mbox{on}\,\gamma,\\
op(B_{2})u=\partial_{n}u_{0}-\partial_{n}v_{0}=e_{2}&\mbox{on}\gamma,
\end{array}
\right.
\end{equation}
 We will proceed like Bellassoued in
\cite{MB}, we will reduce the problem of transmission as a
particular case of a diagonal
 system define only on one side of the interface with boundary conditions.
\par We define the Sobolev spaces with a parameter $\mu$, $H_{\mu}^{s}$ by
$$ u(x,\mu)\in H_{\mu}^{s}\quad \Longleftrightarrow \quad\langle \xi,\mu\rangle^{s} \widehat{u}(\xi,\mu)\in L^{2},\quad
\langle \xi,\mu\rangle^{2}=\abs{\xi}^{2}+\mu^{2},$$
$ \widehat{u}$ denoted the partial Fourier transform with respect to $x$.\\
 For a differential operator
$$ P(x,D,\mu)=\sum_{\abs{\alpha}+k\leq m} a_{\alpha,k}(x)\mu^{k}D^{\alpha},$$
we note the associated symbol by
$$p(x,\xi,\mu)=\sum_{\abs{\alpha}+k\leq m} a_{\alpha,k}(x)\mu^{k}\xi^{\alpha}.$$
The class of symbols of order $m$ is defined by
$$ S_{\mu}^{m}=\left\{ p(x,\xi,\mu)\in C^{\infty}, \abs{D_{x}^{\alpha}D_{\xi}^{\beta}p(x,\xi,\mu)}\leq C_{\alpha,\beta} \langle \xi,\mu\rangle^{m-\abs{\beta}}\right\}$$
and the class of tangential symbols of order $m$ by
$$ \mathcal{T}S_{\mu}^{m}=\left\{ p(x,\xi',\mu)\in C^{\infty}, \abs{D_{x}^{\alpha}D_{\xi'}^{\beta}p(x,\xi',\mu)}\leq C_{\alpha,\beta} \langle \xi',\mu\rangle^{m-\abs{\beta}}\right\}.$$
We denote by $ \mathcal{O}^{m}$ (resp. $ \mathcal{T}\mathcal{O}^{m}$) the set of differentials operators $ P=op (p)$, $ p\in S_{\mu}^{m}$ (resp. $\mathcal{T}S_{\mu}^{m}$).\\
We shall frequently use the symbol $\Lambda=\langle \xi',\mu\rangle=(\abs{\xi'}^{2}+\mu^{2})^{\frac{1}{2}}$.\\
 We shall need
to use the following G{\aa}rding estimate: if $p \in
\mathcal{T}S_{\mu}^{2} $ satisfies for $C_{0}>0$,
$p(x,\xi',\mu)+\overline{p}(x,\xi',\mu)\geq C_{0}\Lambda^{2}$, then
\begin{equation}\label{e3}
\exists\, C_{1}>0,\,\, \exists\, \mu_{0}>0,\,\, \forall \mu
>\mu_{0}, \forall u\in C_{0}^{\infty}(K),\,\, Re(P(x,D',\mu)u,u)\geq
C_{1} \norm{op(\Lambda)u}_{L^{2}}^{2}.
\end{equation}
\par Let $x=(x',x_{n})\in \R^{n-1}\times\R$. In the normal geodesic
system given locally by
$$\Omega_{2}=\{x\in \R^{n}, x_{n}>0\},\quad\quad x_{n}=dist(x, \partial \Omega_{1})=dist(x,x'), $$
the Laplacian is written in the form$$\triangle=- A_{2}(x,D)=-
\left(D_{x_{n}}^{2}+ R(+x_{n},x',D_{x'})\right).$$ The Laplacian on
$\Omega_{1}$ can be identified locally to an operator in
$\Omega_{2}$ gives by $$\triangle= -A_{1}(x,D)=-\left(D_{x_{n}}^{2}+
R(-x_{n},x',D_{x'})\right).$$ We denote the operator, with
$C^{\infty}$ coefficients defined in $\Omega_{2}=\{x_{n}>0\}$, by
$$A(x,D)= \mbox{diag} \Big( A_{1}(x,D_{x}),A_{2}(x,D_{x})\Big)$$
and the tangential operator by
$$ R(x,D_{x'})=\mbox{diag}\Big(R(-x_{n},x',D_{x'}), R(+x_{n},x',D_{x'})\Big)=\mbox{diag}\Big(R_{1}(x,D_{x'}), R_{2}(x,D_{x'})\Big). $$
The principal symbol of the differential operator $A(x,D)$ satisfies
\\ $a(x,\xi)=\xi_{n}^{2}+r(x,\xi')$, where $
r(x,\xi')=\mbox{diag}\Big(r_{1}(x,\xi'),r_{2}(x,\xi') \Big)$ is the
principal symbol of $R(x,D_{x'})$ and the quadratic form
$r_{j}(x,\xi')$, $j=1,2$, satisfies
$$ \exists\, C>0,\quad \forall (x,\xi'),\quad r_{j}(x,\xi')\geq
C \abs{\xi'}^{2},\quad j=1,2.$$ We denote $P(x,D)$ the matrix
operator with $C^{\infty}$ coefficients defined in \\
$\Omega_{2}=\{x_{n}>0\}$, by $$
P(x,D)=\mbox{diag}(P_{1}(x,D),P_{2}(x,D))= \left( \begin{array}{cc}
A_{1}(x,D)-\mu&0\\
0&A_{2}(x,D)-\mu^{2}\\
                                               \end{array}
 \right).
$$
\par Let $\varphi(x)=\mbox{diag}(\varphi_{1}(x),\varphi_{2}(x))$, with $ \varphi_{j}$, $j=1,2$,
are $C^{\infty}$ functions in $ \Omega_{j}$.  For $\mu$ large
enough, we define the operator $$A(x,D,\mu)= e^{\mu
\varphi}A(x,D)e^{-\mu \varphi}:= op(a)$$ where $a\in S_{\mu}^{2}$ is
the principal symbol given by
$$ a(x,\xi,\mu)=\Big(\xi_{n}+i\mu\frac{\partial \varphi}{\partial x_{n}}\Big)^{2}+r\Big(x,\xi'+i\mu\frac{\partial\varphi}{\partial x'}\Big).$$
Let
$$
 op(\tilde{q}_{2,j})=\frac{1}{2}(A_{j}+A_{j}^{\ast}),\quad\quad
op(\tilde{q}_{1,j})=\frac{1}{2i}(A_{j}-A_{j}^{\ast}),\quad\quad
j=1,2 $$ its real and imaginary part. Then we have
\begin{equation}
\left \{\begin{array}{l}\label{e4}
A_{j}=op(\tilde{q}_{2,j})+i op(\tilde{q}_{1,j}),\\\\
\tilde{q}_{2,j}=\xi_{n}^{2}+q_{2,j}(x,\xi',\mu),\quad
\tilde{q}_{1,j}=2\mu\frac{\partial \varphi_{j}}{\partial
x_{n}}\xi_{n}+2\mu q_{1,j}(x,\xi',\mu),\quad j=1,2,
        \end{array}
\right.
\end{equation}
where $q_{1,j}\in \mathcal{T}S_{\mu}^{1}$ and $ q_{2,j}\in
\mathcal{T}S_{\mu}^{2}$ are two tangential symbols given by
\begin{equation}
\left \{\begin{array}{l}\label{e5}
q_{2,j}(x,\xi',\mu)=r_{j}(x,\xi')-(\mu \frac{\partial \varphi_{j}}{\partial x_{n}})^{2}-\mu^{2}r_{j}(x,\frac{\partial \varphi_{j}}{\partial x'}),\\\\
q_{1,j}(x,\xi',\mu)=\tilde{r}_{j}(x,\xi',\frac{\partial
\varphi_{j}}{\partial x'}),\quad j=1,2,
        \end{array}
\right.
\end{equation}
with $\tilde {r}(x,\xi',\eta')$ is the bilinear form associated to
the quadratic form $r(x,\xi')$.\par In the next, $P(x,D,\mu)$ is the
matrix operator with $C^{\infty}$ coefficients defined in
$\Omega_{2}=\{x_{n}>0\}$ by
\begin{equation}\label{e6} P(x,D,\mu)=\mbox{diag}(P_{1}(x,D,\mu),P_{2}(x,D,\mu))= \left( \begin{array}{cc}
A_{1}(x,D,\mu)-\mu&0\\
0&A_{2}(x,D,\mu)-\mu^{2}\\
                                               \end{array}
 \right)
\end{equation}
and $ u=(u_{0},v_{0})$ satisfies the equation
\begin{equation}
\left \{\begin{array}{lc}\label{e16}
Pu=f&\mbox{in}\, \left\{x_{n}>0\right\},\\
op(b_{1})u=u_{0}|_{x_{n}=0}-i\mu v_{0}|_{x_{n}=0}=e_{1}&\mbox{on}\, \left\{x_{n}=0\right\},\\
op(b_{2})u=\left(D_{x_{n}}+i\mu \frac{\partial \varphi_{1}}{\partial
x_{n}}\right)u_{0}|_{x_{n}=0}+\left(D_{x_{n}}+i\mu \frac{\partial
\varphi_{2}}{\partial
x_{n}}\right)v_{0}|_{x_{n}=0}=e_{2}&\mbox{on}\,\left\{x_{n}=0\right\},
\end{array}
\right.
\end{equation}
where $f=(f_{1},f_{2})$, $ e=(e_{1},e_{2})$ and $B=\left(op(b_{1}),
op(b_{2})\right)$. We note $p_{j}(x,\xi,\mu)$, $j=1,2$, the
associated symbol of $P_{j}(x,D,\mu)$.\par We suppose that $\varphi$
satisfies
\begin{equation} \left \{\begin{array}{ll}\label{h1}
\varphi_{1}(x)=\varphi_{2}(x)&
\mbox{on}\,\{x_{n}=0\}\\\\
\dfrac{\partial \varphi_{1}}{\partial x_{n}}>0 &
\mbox{on}\,\{x_{n}=0\}\\\\
\left(\dfrac{\partial \varphi_{1}}{\partial x_{n}}\right)^{2}-
\left(\dfrac{\partial \varphi_{2}}{\partial x_{n}}\right)^{2}>1&
\mbox{on}\,\{x_{n}=0\}
\end{array}
\right.
\end{equation}
and the following condition of hypoellipticity of H\"ormander:
$\exists\, C>0,\,\,\forall x\in K \\\forall\xi\in
\R^{n}\backslash\{0\}, $
\begin{equation}\label{h2}
\,\, \left( \mbox{Re} p_{j}=0\quad\mbox{et}\quad
\frac{1}{2\mu}\mbox{Im}p_{j}=0\right) \,\,\Rightarrow\,\,
\left\{\mbox{Re} p_{j}, \frac{1}{2\mu}\mbox{Im}p_{j}\right\}\geq C
\langle \xi,\mu\rangle^{2},
\end{equation}
where $\{f,g\}(x,\xi)=\sum \left( \frac{\partial f}{\partial
\xi_{j}} \frac{\partial g}{\partial x_{j}} -\frac{\partial
f}{\partial x_{j}}\frac{\partial g}{\partial \xi_{j}} \right) $ is
the Poisson bracket of two functions $f(x,\xi)$ and $g(x,\xi)$ and
$K$ is a compact in $\Omega_{2}$.\par We denote by
$$\norm{u}_{L^{2}(\Omega_{2})}=\norm{u},\quad \norm{u}_{k,\mu}^{2}=\sum_{j=0}^{k}\mu^{2\left( k-j\right) }
\norm{u}_{H^{j}(\Omega_{2})}^{2},
 \quad \norm{u}_{k}^{2}=\norm{op(\Lambda^{k})u}^{2},
$$
$$\abs{u}_{k,\mu}^{2}=\norm{u|_{x_{n}=0}}_{k,\mu}^{2},
\quad \abs{u}_{k}^{2}=\abs{u|_{x_{n}=0}}_{k}^{2},\,\, k\in
\R\quad\mbox{and}\quad\abs{u}_{1,0,\mu}^{2}=
\abs{u}_{1}^{2}+\abs{D_{x_{n}}u}^{2}.$$ We are now ready to state
our result.
\begin{theo}\label{t3}
Let $\varphi$ satisfies (\ref{h1}) and (\ref{h2}). Let $w\in C_{0}^{\infty}(\overline{\Omega}_{2})
 $ and $ \chi\in
C_{0}^{\infty}(\R^{n+1})$ such that $\chi=1$ in the support of $w$. Then there exists
constants $C>0$ and $\mu_{0}>0$ such that for any $\mu\geq\mu_{0}$
we have the following estimate
\begin{eqnarray}\label{tt}
\mu\norm{w}_{1,\mu}^{2}+\mu^{2}\abs{w}^{2}_{\frac{1}{2}}+\mu^{2}\abs{D_{x_{n}}w}^{2}_{-\frac{1}{2}}~~~~~~~~~~~~~~~~~~~~
~~~~~~~~~~~~~~~~~~~~~~~~~~~~~~~\nonumber\\\nonumber\\
\leq C\left( \norm{P(x,D,\mu)w}^{2}+
\abs{op(b_{1})w}^{2}_{\frac{1}{2}}+\mu\abs{op(b_{2})w}^{2}\right).
\end{eqnarray}
\end{theo}
\begin{coro}\label{t4}
Let $\varphi$ satisfies (\ref{h1}) and (\ref{h2}). Then there exists
constants $C>0$ and $\mu_{0}>0$ such that for any $\mu\geq\mu_{0}$
we have the following estimate
\begin{equation}\label{t41}
\mu \norm{e^{\mu \varphi}h}^{2}_{H^{1}}\leq C\left( \norm{e^{\mu
\varphi}P(x,D)h}^{2}+\abs{e^{\mu\varphi}op(B_{1})h}^{2}_{H^{\frac{1}{2}}}+\mu
\abs{e^{\mu\varphi}op(B_{2})h}^{2}\right),
\end{equation}
for any $h\in C_{0}^{\infty}(\overline{\Omega}_{2})$.
\end{coro}
{\bf{Proof.}}\\
Let $w=e^{\mu\varphi}h$. Recalling that $
P(x,D,\mu)w=e^{\mu\varphi}P(x,D)e^{-\mu\varphi}w$ and using
(\ref{tt}), we obtain (\ref{t41}).
\subsection{Proof of Theorem \ref{II}}
\hspace{5mm} We denote $x=(x',x_{n})$ a point in $\Omega$. Let
$x_{0}=(0,-\delta)$, $\delta>0$. We set $$
\psi(x)=\abs{x-x_{0}}^{2}-\delta^{2}\qquad \mbox{and}$$
$$ \varphi_{1}(x)=e^{-\beta\psi(x',-x_{n})}, \quad \varphi_{2}(x)=e^{-\beta(\psi(x)-\alpha x_{n})},\quad \beta>0,\quad\mbox{and}\quad
\frac{\delta}{2}<\alpha<2\delta.$$ The weight function $
\varphi=diag(\varphi_{1},\varphi_{2})$ has to satisfy (\ref{h1}) and
(\ref{h2}). With these choices, we have
$\varphi_{1}|_{x_{n}=0}=\varphi_{2}|_{x_{n}=0}$ and $\frac{\partial
\varphi_{1}}{\partial x_{n}}|_{x_{n}=0}>0$. It remains to verify
\begin{equation}\label{1c}
\left(\frac{\partial \varphi_{1}}{\partial x_{n}} \right)^{2}-
   \left(\frac{\partial \varphi_{2}}{\partial x_{n}}
   \right)^{2}>1\quad \mbox{on}\,\{x_{n}=0\}
\end{equation}
and the condition (\ref{h2}).
 We begin by condition (\ref{h2}) and we compute for $\varphi_{1}$ and $p_{1}$ (the computation for $\varphi_{2}$ and $p_{2}$ is made in the same way).
 Recalling that
\begin{eqnarray*}
\left\{\mbox{Re} p_{1},
\frac{1}{2\mu}\mbox{Im}p_{1}\right\}(x,\xi)=\frac{\mbox{Im}}{2\mu}\left[
\partial_{\xi}p_{1}(x,\xi-i\mu \varphi'_{1}(x))\,\partial_{x}p_{1}(x,\xi+i\mu \varphi'_{1}(x))\right]
\\\\+^{t}\left[ \partial_{\xi}p_{1}(x,\xi-i\mu
\varphi'_{1}(x))\right]\varphi''_{1}(x)\left[\partial_{\xi}p_{1}(x,\xi-i\mu
\varphi'_{1}(x)) \right].
\end{eqnarray*}
We replace $\varphi_{1}(x)$ by $\varphi_{1}(x)= e^{-\beta
\psi(x',-x_{n})}$, $\beta>0$, we obtain, by noting $\xi=-\beta
\varphi_{1}(x)\eta$
\begin{eqnarray*}
\left\{\mbox{Re} p_{1},
\frac{1}{2\mu}\mbox{Im}p_{1}\right\}(x,\xi)~~~~~~~~~~~~~~~~~~~
~~~~~~~~~~~~~~~~~~~~~~~~~~~~~~~~~~~~~~~~~~~~~~~~~~~~~~~~~~~~~~~~~~~~~~~~\\\\
=(-\beta \varphi_{1})^{3}\left[\left\{\mbox{Re}
p_{1}(x,\eta-i\mu\psi'),\frac{1}{2\mu}\mbox{Im}p_{1}(x,\eta+i\mu\psi')\right\}(x,\eta)
  -\beta \left|\psi'(x) \partial_{\eta}p_{1}(x,\eta+i\mu\psi')\right|^{2}\right]
\end{eqnarray*}
and $$\left|\psi'(x)
\partial_{\eta}p_{1}(x,\eta+i\mu\psi')\right|^{2}=4\left[ \mu^{2}\abs{p_{1}(x,\psi')}^{2}+
\abs{\tilde{p}_{1}(x,\eta,\psi')}^{2}\right]
$$
where $ \tilde{p}_{1}(x,\eta,\psi')$ is the bilinear form associated
to the quadratic form $p_{1}(x,\eta)$. We have
$$\left( \mbox{Re} p_{1}=0\quad\mbox{et}\quad
\frac{1}{2\mu}\mbox{Im}p_{1}=0\right) \Longleftrightarrow
p_{1}(x,\eta+i\mu\psi')=0,
$$
\begin{itemize}
\item If $\mu= 0$, we have $ p_{1}(x,\xi)=0$ which is impossible.
Indeed, we have\\ $ p_{1}(x,\xi)\geq C\abs{\xi}^{2}$, $\forall
(x,\xi)\in K \times\R^{n}$, $K$ compact in $\Omega_{2}$.
\item If $\mu\neq 0$, we have $ \tilde{p}_{1}(x,\eta,\psi')=0$. \\ Then
$ \left|\psi'(x)
\partial_{\eta}p_{1}(x,\eta+i\mu\psi')\right|^{2}=4\mu^{2}\abs{p_{1}(x,\psi')}^{2}>0$.
On the other hand, we have $$ \left\{\mbox{Re}
p_{1}(x,\eta-i\mu\psi'),\frac{1}{2\mu}\mbox{Im}p_{1}(x,\eta+i\mu\psi')\right\}(x,\eta)\leq
C_{1}(\abs{\eta}^{2}+\mu^{2}\abs{\psi'}^{2})$$ where $C_{1}$ is a
positive constant independent of $\psi'$. Then for $\beta\geq
C_{1}$, we satisfy the condition (\ref{h2}).
\end{itemize}
Now let us verify (\ref{1c}). We have, on $\{x_{n}=0\}$,
$$\left(\frac{\partial \varphi_{1}}{\partial x_{n}} \right)^{2}-
   \left(\frac{\partial \varphi_{2}}{\partial x_{n}}
   \right)^{2}=\beta^{2}\alpha(4\delta-\alpha)e^{-2\beta\psi}. $$
 Then to satisfy (\ref{1c}), it suffices to choose
   $\beta=\frac{M}{\delta}$ where $M>0$ such that $\frac{M}{\delta}\geq
   C_{1}$.
\par We now choose
$r_{1}<r'_{1}<r_{2}<0=\psi(0)<r'_{2}<r_{3}<r'_{3}$. We denote
$$w_{j}=\{ x\in\Omega,\, r_{j}<\psi(x)<r'_{j}\} \quad\mbox{and}\quad T_{x_{0}}=w_{2}\cap\Omega_{2}.$$
We set $R_{j}=e^{-\beta r_{j}}$, $ R'_{j}=e^{-\beta r'_{j}} $,
$j=1,2,3$.\\ Then $ R'_{3}<R_{3}<R'_{2}<R_{2}<R'_{1}<R_{1}$. We need
also to introduce a cut-off function $ \tilde{\chi}\in
C_{0}^{\infty}(\R^{n+1})$ such that
$$\tilde{\chi}(\rho)=\left \{\begin{array}{lcc}
0&\mbox{if}& \rho\leq r_{1},\quad \rho\geq r'_{3}\\\\
1&\mbox{if}& \rho \in [r'_{1}, r_{3}].\end{array} \right. $$ Let
$\tilde{u}=(
\tilde{u}_{0},\tilde{v}_{0})=\tilde{\chi}u=(\tilde{\chi}u_{0},\tilde{\chi}v_{0})$.
Then we get the following system
$$ \left \{\begin{array}{rcl}
(\triangle-i\mu) \tilde{u}_{0}&=&\tilde{\chi}f_{0}+ [\triangle-i\mu,\tilde{\chi}]u_{0}\\
(\triangle+\mu^{2})\tilde{v}_{0}&=&\tilde{\chi}(g_{1}+i\mu g_{0})+ [\triangle+\mu^{2}, \tilde{\chi}]v_{0},\\
\tilde{v}_{1}&=&g_{0}+i\mu \tilde{v}_{0},
\end{array} \right.$$
with the following boundary conditions
$$ \left \{\begin{array}{lcl}
\tilde{u}_{0}|_{\Gamma_{1}}=\tilde{v}_{0}|_{\Gamma_{2}}&=&0,\\
op(b_{1})\tilde{u}=\tilde{u}_{0}-i\mu \tilde{v}_{0}&=&(\tilde{\chi}g_{0})|_{\gamma},\\
op(b_{2})\tilde{u}&=&([\partial_{n},\tilde{\chi}]u_{0}-[\partial_{n},
\tilde{\chi}]v_{0})|_{\gamma}.
\end{array} \right.$$
From the Carleman estimate of Corollary \ref{t4} , we have
\begin{equation}\label{s3}
\mu\norm{e^{\mu\varphi}\tilde{u}}_{H^{1}}^{2}\leq
C\left(\norm{e^{\mu \varphi_{1}}( \triangle-i\mu)\tilde{u}_{0}}^{2}+
\norm{e^{\mu \varphi_{2}}( \triangle+\mu^{2})\tilde{v}_{0}}^{2}+
\abs{e^{\mu
\varphi}op(b_{1})\tilde{u}}_{H^{\frac{1}{2}}}^{2}+\mu\abs{e^{\mu
\varphi}op(b_{2})\tilde{u}}^{2} \right).
\end{equation}
Using the fact $[\triangle-i\mu,\tilde{\chi}]$ is the first order
operator supported in $ (w_{1} \cup w_{3}) \cap \Omega_{1}$, we have
\begin{equation}\label{s4}
\norm{e^{\mu \varphi_{1}}( \triangle-i\mu)\tilde{u}_{0}}^{2}\leq
C\left( e^{2\mu R_{1}}\norm{f_{0}}^{2}_{L^{2}(\Omega_{1})}+ e^{2\mu
R_{1}}\norm{u_{0}}^{2}_{H^{1}(\Omega_{1})}\right).
\end{equation}
Recalling that $ [ \triangle+\mu^{2},\tilde{\chi}]$ is the first
order operator supported in $(w_{1} \cup w_{3})\cap \Omega_{2}$, we
show
\begin{equation}\label{s5}
\norm{e^{\mu \varphi_{2}}( \triangle+\mu^{2})\tilde{v}_{0}}^{2}\leq
C\left(e^{2\mu}\norm{ g_{1}+i\mu
g_{0}}^{2}_{L^{2}(\Omega_{2})}+e^{2\mu
R_{3}}\norm{v_{0}}_{H^{1}(\Omega_{2})}^{2} \right).
\end{equation}
From the trace formula and recalling that $op(b_{2})\tilde{u}$ is an
operator of order zero and supported in $\{x_{n}=0\}\cap w_{3}$, we
show
\begin{equation}\label{s6}
\mu\abs{e^{\mu \varphi}op(b_{2})\tilde{u}}^{2}\leq C e^{2\mu
R_{3}}\norm{u}_{H^{1}(\Omega)}^{2}\leq C\left(e^{2\mu
R_{3}}\norm{u_{0}}_{H^{1}(\Omega_{1})}^{2}+e^{2\mu
R_{3}}\norm{v_{0}}_{H^{1}(\Omega_{2})}^{2} \right).
\end{equation}
Now we need to use this result shown in Appendix B
\begin{lemm}\label{le2}
There exists $C>0$ such that for all $s\in\R$ and $u\in
C_{0}^{\infty}(\Omega)$, we have
\begin{equation}\label{le21}
\norm{op(\Lambda^{s})e^{\mu\varphi}u}\leq C e^{\mu
C}\norm{op(\Lambda^{s})u}.
\end{equation}
\end{lemm} Following this Lemma, we obtain
\begin{equation}\label{s7}
\abs{e^{\mu \varphi}op(b_{1})\tilde{u}}_{H^{\frac{1}{2}}}^{2}\leq C
e^{2\mu c}\abs{g_{0}}_{H^{\frac{1}{2}}}^{2}\leq C e^{2\mu
c}\norm{g_{0}}^{2}_{H^{1}(\Omega_{2})}.
\end{equation}
Combining  (\ref{s3}), (\ref{s4}), (\ref{s5}), (\ref{s6}) and
(\ref{s7}), we obtain
\begin{eqnarray*}
C\mu e^{2\mu R'_{2}}\norm{u_{0}}^{2}_{H^{1}(w_{2}\cap\Omega_{1})}+
C\mu e^{2\mu R'_{2}}\norm{v_{0}}^{2}_{H^{1}(T_{x_{0}})}\leq
C(e^{2\mu R_{1}}\norm{f_{0}}^{2}_{L^{2}(\Omega_{1})}+ e^{2\mu
R_{1}}\norm{u_{0}}^{2}_{H^{1}(\Omega_{1})}~~\\\\
+e^{2\mu}\norm{ g_{1}+i\mu g_{0}}^{2}_{L^{2}(\Omega_{2})}+e^{2\mu
R_{3}}\norm{v_{0}}_{H^{1}(\Omega_{2})}^{2}+e^{2\mu
R_{3}}\norm{u_{0}}_{H^{1}(\Omega_{1})}^{2}+e^{2\mu
c}\norm{g_{0}}^{2}_{H^{1}(\Omega_{2})}).
\end{eqnarray*}
Since $R'_{2}<R_{1}$. Then we have
\begin{eqnarray}\label{s8}
\norm{v_{0}}^{2}_{H^{1}(T_{x_{0}})}&\leq& C
e^{c_{1}\mu}\left[\norm{f_{0}}^{2}_{L^{2}(\Omega_{1})}+
 \norm{ g_{1}+i\mu
g_{0}}^{2}_{L^{2}(\Omega_{2})}+\norm{g_{0}}^{2}_{H^{1}(\Omega_{2})}+
\norm{u_{0}}^{2}_{H^{1}(\Omega_{1})}\right]\nonumber\\\nonumber\\&+&C
e^{-c_{2}\mu}\norm{v_{0}}_{H^{1}(\Omega_{2})}^{2}.
\end{eqnarray}
Since $\gamma$ is compact, then there exists a finite number of
$T_{x_{0}}$. Let $V_{\epsilon_{2}} \subset\cup T_{x_{0}}$. Then we
obtain (\ref{II1})
\subsection{Proof of Carleman estimate (Theorem \ref{t3})}
\hspace{5mm}In the first step, we state the following estimates
\begin{theo}\label{t2}
Let $ \varphi$ satisfies (\ref{h1}) and (\ref{h2}). Then there
exists constants $C>0$ and $\mu_{0}$ such that for any
$\mu\geq\mu_{0}$ we have the following estimates
\begin{equation}\label{k4}
\mu\norm{u}_{1,\mu}^{2}\leq C\left( \norm{P(x,D,\mu)u}^{2}+\mu
\abs{u}_{1,0,\mu}^{2}\right)
\end{equation}
and
\begin{equation}\label{tt2}
\mu\norm{u}^{2}_{1,\mu}+\mu\abs{u}_{1,0,\mu}^{2}\leq C\left(
\norm{P(x,D,\mu)u}^{2}+\mu^{-1}\abs{op(b_{1})u}^{2}_{1}+\mu\abs{op(b_{2})u}^{2}\right),
\end{equation}
for any $u\in C_{0}^{\infty}(\overline{\Omega}_{2})$.
\end{theo}
In the second step, we need to prove this Lemma
\begin{lemm}\label{ll}
There exists constants $C>0$ and $\mu_{0}>0$ such that for any
$\mu\geq\mu_{0}$ we have the following estimate
\begin{eqnarray}\label{lll}
\norm{D_{x_{n}}^{2}op( \Lambda^{-\frac{1}{2}})u}^{2}+
\norm{D_{x_{n}}op(
\Lambda^{\frac{1}{2}})u}^{2}+\norm{op(\Lambda^{\frac{3}{2}})u}^{2}+\mu\abs{u}^{2}_{1,0,\mu}
~~~~\nonumber\\\nonumber\\
\leq C\left( \norm{
P(x,D,\mu)u}^{2}+\mu^{-1}\abs{op(b_{1})u}^{2}_{1}+\mu\abs{op(b_{2})u}^{2}\right),
\end{eqnarray}
for any  $ u\in C_{0}^{\infty}(\overline{\Omega}_{2})$.
\end{lemm}
{\bf Proof.}\\
We have $$P(x,D,\mu)=D_{x_{n}}^{2}+R+\mu C_{1}+\mu^{2}C_{0},$$ where $R\in
\mathcal{T}\mathcal{O}^{2}$, $C_{1}=c_{1}(x)D_{x_{n}}+T_{1}$, with
$T_{1}\in \mathcal{T}\mathcal{O}^{1}$ and $C_{0}\in
\mathcal{T}\mathcal{O}^{0}$. Then we have
\begin{eqnarray*}\norm{(D_{x_{n}}^{2}+R)op(\Lambda^{-\frac{1}{2}})u}^{2}~~~~~~~~~~~~~~~~~~~~~~~~~~~~~~~~~~
~~~~~~~~~~~~~~~~~~~~~~~~~~~~~~~~~~~~~~~~~~~~~~~~\\\leq C\left(
\norm{Pop(\Lambda^{-\frac{1}{2}})u}^{2}
+\mu^{2}\norm{op(\Lambda^{\frac{1}{2}})u}^{2}+\mu^{2}\norm{D_{x_{n}}op(\Lambda^{-\frac{1}{2}})u}^{2}+\mu^{4}
\norm{op(\Lambda^{-\frac{1}{2}})u}^{2}\right).
\end{eqnarray*}
Since $$\begin{array} {l}\mu^{4}
\norm{op(\Lambda^{-\frac{1}{2}})u}^{2}\leq C\mu^{3}\norm{u}^{2},\\
\mu^{2}\norm{D_{x_{n}}op(\Lambda^{-\frac{1}{2}})u}^{2}\leq C
\mu\norm{D_{x_{n}}u}^{2}\quad\mbox{and}\\
\mu^{2}\norm{op(\Lambda^{\frac{1}{2}})u}^{2}=\mu^{2}(\frac{1}{\sqrt{\mu}}op(\Lambda)u,\sqrt{\mu}u)\leq
C \left( \mu
\norm{op(\Lambda)u}^{2}+\mu^{3}\norm{u}^{2}\right).\end{array}$$
Using the fact that $\norm{u}^{2}_{1,\mu}\simeq
\norm{op(\Lambda)u}^{2}+\norm{D_{x_{n}}u}^{2}$, we obtain
$$\norm{(D_{x_{n}}^{2}+R)op(\Lambda^{-\frac{1}{2}})u}^{2}\leq C\left(
\norm{Pop(\Lambda^{-\frac{1}{2}})u}^{2}+\mu\norm{u}_{1,\mu}^{2}\right).$$
Following (\ref{k4}), we have
\begin{equation}\label{ll1}
\norm{(D_{x_{n}}^{2}+R)op(\Lambda^{-\frac{1}{2}})u}^{2}\leq C\left(
\norm{Pop(\Lambda^{-\frac{1}{2}})u}^{2}+\norm{Pu}^{2}+\mu\abs{u}^{2}_{1,0,\mu}\right).
\end{equation}
We can write
\begin{eqnarray}\label{ll2}
  Pop(\Lambda^{-\frac{1}{2}})u &=& op(\Lambda^{-\frac{1}{2}})Pu+[P, op(\Lambda^{-\frac{1}{2}})]u \nonumber \\
   &=& op(\Lambda^{-\frac{1}{2}})Pu+[R, op(\Lambda^{-\frac{1}{2}})]u \nonumber\\
   &+&  \mu [C_{1}, op(\Lambda^{-\frac{1}{2}})]u
   + \mu^{2}[C_{0}, op(\Lambda^{-\frac{1}{2}})]u \nonumber\\
   &=&op(\Lambda^{-\frac{1}{2}})Pu+t_{1}+t_{2}+t_{3}.
\end{eqnarray}
Let us estimate $t_{1}$, $t_{2}$ and $t_{3}$. We have $[R,
op(\Lambda^{-\frac{1}{2}})]\in \mathcal{T}\mathcal{O}^{\frac{1}{2}}
$, then following (\ref{k4}), we have
\begin{equation}\label{ll4} \norm{t_{1}}^{2}\leq C\norm{
op(\Lambda^{\frac{1}{2}})u}^{2}\leq C\left(\norm{
op(\Lambda)u}^{2}+\norm{u}^{2}\right) \leq C\left(
\norm{Pu}^{2}+\mu\abs{u}^{2}_{1,0,\mu}\right).
\end{equation}
We have $t_{3}=\mu [C_{1}, op(\Lambda^{-\frac{1}{2}})]u=\mu
[c_{1}(x)D_{x_{n}}, op(\Lambda^{-\frac{1}{2}})]u+ \mu [T_{1},
op(\Lambda^{-\frac{1}{2}})]u$. Then following (\ref{k4}), we obtain
\begin{equation}\label{ll5} \norm{t_{2}}^{2}\leq C\left(
\mu^{-1}\norm{D_{x_{n}}u}^{2}+\mu\norm{u}^{2}\right)\leq
C\left(\norm{Pu}^{2}+\mu\abs{u}^{2}_{1,0,\mu}\right).
\end{equation}
We have $ [C_{0}, op(\Lambda^{-\frac{1}{2}})]\in
\mathcal{T}\mathcal{O}^{-\frac{3}{2}}$, then following (\ref{k4}),
we obtain \begin{equation}\label{ll6} \norm{ \mu^{2}[C_{0},
op(\Lambda^{-\frac{1}{2}})]u}^{2}\leq C\mu \norm{u}^{2}\leq
C\left(\norm{Pu}^{2}+\mu\abs{u}^{2}_{1,0,\mu}\right)
\end{equation}
From (\ref{ll2}), (\ref{ll4}), (\ref{ll5}) and (\ref{ll6}), we have
$$\norm{  Pop(\Lambda^{-\frac{1}{2}})u}^{2}\leq C\left(\norm{Pu}^{2}+\mu\abs{u}^{2}_{1,0,\mu}\right). $$
Inserting this inequality in (\ref{ll1}), we obtain
\begin{equation}\label{ll7}
\norm{(D_{x_{n}}^{2}+R)op(\Lambda^{-\frac{1}{2}})u}^{2}\leq
C\left(\norm{Pu}^{2}+\mu\abs{u}^{2}_{1,0,\mu}\right).
\end{equation}
Moreover, we have
$$\norm{(D_{x_{n}}^{2}+R)op(\Lambda^{-\frac{1}{2}})u}^{2}=\norm{D_{x_{n}}^{2}op(\Lambda^{-\frac{1}{2}})u}^{2}+
\norm{Rop(\Lambda^{-\frac{1}{2}})u}^{2}+2\mathcal{R}e(D_{x_{n}}^{2}op(\Lambda^{-\frac{1}{2}})u,Rop(\Lambda^{-\frac{1}{2}})u),
$$
where $(., .)$ denoted the scalar product in $L^{2}$.
By integration by parts, we find
\begin{eqnarray}\label{ll10}
\norm{(D_{x_{n}}^{2}+R)op(\Lambda^{-\frac{1}{2}})u}^{2}=\norm{D_{x_{n}}^{2}op(\Lambda^{-\frac{1}{2}})u}^{2}+
\norm{Rop(\Lambda^{-\frac{1}{2}})u}^{2}~~~~~~~~~~~~~~~~~~~~~~~~~~~~~~~\nonumber\\\nonumber\\+2\mathcal{R}e\left(i(D_{x_{n}}u,R
op(\Lambda^{-1})u)_{0}+i(D_{x_{n}}u, [ op(\Lambda^{-\frac{1}{2}}),R]
op(\Lambda^{-\frac{1}{2}})u)_{0}
\right)~~~~~~~~~~~~~~~~~~~~\nonumber\\\nonumber\\+2\mathcal{R}e\left((RD_{x_{n}}op(\Lambda^{-\frac{1}{2}})u,D_{x_{n}}op(\Lambda^{-\frac{1}{2}})u
)+ (
D_{x_{n}}op(\Lambda^{-\frac{1}{2}})u,[D_{x_{n}},R]op(\Lambda^{-\frac{1}{2}})u)\right)
.
\end{eqnarray}
Since, we have
$$ \norm{op(\Lambda^{\frac{3}{2}})u}^{2}= (op(\Lambda^{2})op(\Lambda^{\frac{1}{2}})u, op(\Lambda^{\frac{1}{2}})u)=
\sum_{j\leq n-1}(D_{j}^{2}op(\Lambda^{\frac{1}{2}})u,
op(\Lambda^{\frac{1}{2}})u )+ \mu^{2} ( op(\Lambda^{\frac{1}{2}})u,
op(\Lambda^{\frac{1}{2}})u).$$ By integration by parts, we find
\begin{equation}\label{ll14}
\norm{op(\Lambda^{\frac{3}{2}})u}^{2}= \sum_{j\leq
n-1}(D_{j}op(\Lambda^{\frac{1}{2}})u,D_{j}
op(\Lambda^{\frac{1}{2}})u )+ \mu^{2} \norm{
op(\Lambda^{\frac{1}{2}})u}^{2}= k+\mu^{2} \norm{
op(\Lambda^{\frac{1}{2}})u}^{2} .
\end{equation}
Let $\chi_{0}\in C_{0}^{\infty}(\R^{n+1})$ such that $ \chi_{0}=1$
in the support of $u$. We have $$k=\sum_{j\leq
n-1}(\chi_{0}D_{j}op(\Lambda^{\frac{1}{2}})u,D_{j}
op(\Lambda^{\frac{1}{2}})u )+ \sum_{j\leq
n-1}((1-\chi_{0})D_{j}op(\Lambda^{\frac{1}{2}})u,D_{j}
op(\Lambda^{\frac{1}{2}})u ).$$ Recalling that $\chi_{0}u=u$, we
obtain
\begin{equation}\label{ll15}
k=\sum_{j\leq n-1}(\chi_{0}D_{j}op(\Lambda^{\frac{1}{2}})u,D_{j}
op(\Lambda^{\frac{1}{2}})u )+\sum_{j\leq
n-1}([(1-\chi_{0}),D_{j}op(\Lambda^{\frac{1}{2}})]u,D_{j}
op(\Lambda^{\frac{1}{2}})u )=k'+k".
\end{equation}
Using the fact that $
[(1-\chi_{0}),D_{j}op(\Lambda^{\frac{1}{2}})]\in
\mathcal{T}\mathcal{O}^{\frac{1}{2}}$ and $ D_{j}
op(\Lambda^{\frac{1}{2}})\in \mathcal{T}\mathcal{O}^{\frac{3}{2}}$,
we show
\begin{equation}\label{ll16}
k"\leq C\norm{op(\Lambda)u }^{2}.
\end{equation}
Using the fact that $\sum_{j,k\leq n-1}\chi_{0}a_{j,k}D_{j}v
\overline{D_{k}v }\geq \delta \chi_{0}\sum_{j\leq
n-1}\abs{D_{j}v}^{2}$, $\delta >0$, we obtain
\begin{eqnarray*}
k'\leq C\sum_{j,k\leq n-1}(\chi_{0}
a_{jk}D_{j}op(\Lambda^{\frac{1}{2}})u,D_{k}
op(\Lambda^{\frac{1}{2}})u
)~~~~~~~~~~~~~~~~~~~~~~~~~~~~~~~~~~~~~~~~~~~~~~~~~~~~~~~~~
\\\\
\leq C\sum_{j,k\leq n-1}([\chi_{0},
a_{jk}D_{j}op(\Lambda^{\frac{1}{2}})]u,D_{k}
op(\Lambda^{\frac{1}{2}})u )+\sum_{j,k\leq
n-1}(a_{jk}D_{j}op(\Lambda^{\frac{1}{2}})u,D_{k}
op(\Lambda^{\frac{1}{2}})u ).
\end{eqnarray*}
Using the fact that $[\chi_{0},
a_{jk}D_{j}op(\Lambda^{\frac{1}{2}})]\in
\mathcal{T}\mathcal{O}^{\frac{1}{2}} $ and $ D_{k}
op(\Lambda^{\frac{1}{2}})u \in
\mathcal{T}\mathcal{O}^{\frac{3}{2}}$, we obtain
\begin{equation}\label{ll17}
k'\leq C\left(\sum_{j,k\leq
n-1}(a_{jk}D_{j}op(\Lambda^{\frac{1}{2}})u,D_{k}
op(\Lambda^{\frac{1}{2}})u )+\norm{op(\Lambda)u }^{2}\right).
\end{equation}
By integratin by parts and recalling that $R=\sum_{j,k\leq
n-1}a_{j,k}D_{j}D_{k}$, we have
\begin{eqnarray}\label{ll18}
\sum_{j,k\leq n-1}(a_{jk}D_{j}op(\Lambda^{\frac{1}{2}})u,D_{k}
op(\Lambda^{\frac{1}{2}})u )&=& (Rop(\Lambda^{\frac{1}{2}})u,
op(\Lambda^{\frac{1}{2}})u \\&+&\sum_{j,k\leq n-1}([ D_{k}, a_{jk}]
D_{j}op(\Lambda^{\frac{1}{2}})u,
op(\Lambda^{\frac{1}{2}})u)\nonumber.
\end{eqnarray}
Since $[ D_{k}, a_{jk}] D_{j}op(\Lambda^{\frac{1}{2}})\in
\mathcal{T}\mathcal{O}^{\frac{3}{2}} $, then
$$ \sum_{j,k\leq n-1}([ D_{k}, a_{jk}]
D_{j}op(\Lambda^{\frac{1}{2}})u, op(\Lambda^{\frac{1}{2}})u)\leq
C\norm{op(\Lambda)u }^{2}.
$$
Following (\ref{ll18}), we obtain
\begin{equation}\label{ll19}
\sum_{j,k\leq n-1}(a_{jk}D_{j}op(\Lambda^{\frac{1}{2}})u,D_{k}
op(\Lambda^{\frac{1}{2}})u )\leq
C\left((Rop(\Lambda^{\frac{1}{2}})u, op(\Lambda^{\frac{1}{2}})u)
+\norm{op(\Lambda)u }^{2}\right).
\end{equation}
Since
$$ (Rop(\Lambda^{\frac{1}{2}})u, op(\Lambda^{\frac{1}{2}})u)=
(Rop(\Lambda^{-\frac{1}{2}})u, op(\Lambda^{\frac{3}{2}})u)+ (
[op(\Lambda^{-1}), R]op(\Lambda^{\frac{1}{2}})u,
op(\Lambda^{\frac{3}{2}})u) ).$$ Using the fact that $
[op(\Lambda^{-1}), R]op(\Lambda^{\frac{1}{2}})\in
\mathcal{T}\mathcal{O}^{\frac{1}{2}}$ and the Cauchy Schwartz
inequality, we obtain
\begin{equation}\label{ll20}
(Rop(\Lambda^{\frac{1}{2}})u, op(\Lambda^{\frac{1}{2}})u)\leq
\epsilon C\norm{op(\Lambda^{\frac{3}{2}})u
}^{2}+\frac{C}{\epsilon}\norm{Rop(\Lambda^{-\frac{1}{2}})u}^{2}
+C\norm{op(\Lambda)u}^{2}
\end{equation}
Combining (\ref{ll14}), (\ref{ll15}), (\ref{ll16}), (\ref{ll17}),
(\ref{ll19}) and (\ref{ll20}), we obtain
$$ \norm{op(\Lambda^{\frac{3}{2}})u
}^{2}\leq \epsilon C\norm{op(\Lambda^{\frac{3}{2}})u
}^{2}+\frac{C}{\epsilon}\norm{Rop(\Lambda^{-\frac{1}{2}})u}^{2}
+C\norm{op(\Lambda)u}^{2}.$$ For $ \epsilon$ small enough, we obtain
\begin{equation}\label{ll11}
\norm{Rop(\Lambda^{-\frac{1}{2}})u}^{2}\geq
C\left(\norm{op(\Lambda^{\frac{3}{2}})u }^{2}-
\mu^{2}\norm{op(\Lambda^{\frac{1}{2}})u}^{2} \right).
\end{equation}
Using the same computations, we show
\begin{equation}\label{ll12}
(RD_{x_{n}}op(\Lambda^{-\frac{1}{2}})u,D_{x_{n}}op(\Lambda^{-\frac{1}{2}})u)\geq
C\left( \norm{D_{x_{n}}op(\Lambda^{\frac{1}{2}})u }^{2}-\mu\norm{
D_{x_{n}}u}^{2}\right).
\end{equation}
Combining (\ref{ll10}), (\ref{ll11}) and (\ref{ll12}), we obtain
\begin{eqnarray}\label{ll8}
\norm{(D_{x_{n}}^{2}+R)op(\Lambda^{-\frac{1}{2}})u}^{2}&+&\abs{(D_{x_{n}}u,R
op(\Lambda^{-1})u)_{0}}+ \abs{(D_{x_{n}}u, [
op(\Lambda^{-\frac{1}{2}}),R]
op(\Lambda^{-\frac{1}{2}})u)_{0}}\nonumber\\\nonumber\\
&+& \abs{(
D_{x_{n}}op(\Lambda^{-\frac{1}{2}})u,[D_{x_{n}},R]op(\Lambda^{-\frac{1}{2}})u)}
+\mu\norm{u}_{1,\mu}^{2}\\\nonumber\\
&\geq& C\left( \norm{D_{x_{n}}^{2}op(\Lambda^{-\frac{1}{2}})u}^{2}+
\norm{D_{x_{n}}
op(\Lambda^{\frac{1}{2}})u}^{2}+\norm{op(\Lambda^{\frac{3}{2}})u}^{2}\right)\nonumber.
\end{eqnarray}
Since
\begin{equation}\label{ll9}
\abs{(D_{x_{n}}u,R op(\Lambda^{-1})u)_{0}}+ \abs{(D_{x_{n}}u, [
op(\Lambda^{-\frac{1}{2}}),R] op(\Lambda^{-\frac{1}{2}})u)_{0}}
 \leq C\left(
\abs{D_{x_{n}}u}^{2}+\abs{u}_{1}^{2}\right)=C\abs{u}_{1,0,\mu}^{2}
\end{equation}
and
\begin{equation}\label{ll13}
\abs{(
D_{x_{n}}op(\Lambda^{-\frac{1}{2}})u,[D_{x_{n}},R]op(\Lambda^{-\frac{1}{2}})u)}\leq
C\mu\norm{u}_{1,\mu}^{2}.
\end{equation}
 From (\ref{ll7}), (\ref{ll8}), (\ref{ll9}), (\ref{ll13}) and
 (\ref{k4}),  we obtain
\begin{eqnarray*}
\norm{D_{x_{n}}^{2}op(\Lambda^{-\frac{1}{2}})u}^{2}+ \norm{D_{x_{n}}
op(\Lambda^{\frac{1}{2}})u}^{2}+\norm{op(\Lambda^{\frac{3}{2}})u}^{2}~~~~~~~~~~~~~~\\\\
\leq C\left( \norm{P(x,D,\mu)u}^{2}+\mu\abs{u}_{1,0,\mu}^{2}\right).
\end{eqnarray*}
Following (\ref{tt2}), we obtain (\ref{lll}).\\
$~~~~~~~~~~~~~~~~~~~~~~~~~~~~~~~~~~~~~~~~~~~~~~~~~~~~~~~~~~~~~~
~~~~~~~~~~~~~~~~~~~~~~~~~~~~~~~~~~~~~~~~~~~~~~~~~~~~~~~~~~~~\square$
\par We are now ready to prove Theorem \ref{t3}.\par Let $ \chi\in C_{0}^{\infty}(\R^{n+1})$ such that $\chi=1$ in
the support of $w$ and $u=\chi op(\Lambda^{-\frac{1}{2}})w$. Then
\begin{eqnarray}\label{ttt1}
Pu&=&op(\Lambda^{-\frac{1}{2}})Pw+
[P,op(\Lambda^{-\frac{1}{2}})]w+P[\chi,op(\Lambda^{-\frac{1}{2}})]w
\nonumber\\
&=&op(\Lambda^{-\frac{1}{2}})Pw +[P,op(\Lambda^{-\frac{1}{2}})]w
+D_{x_{n}}^{2}[\chi,op(\Lambda^{-\frac{1}{2}})]w \nonumber\\
&+& R[\chi,op(\Lambda^{-\frac{1}{2}})]w + \mu
c_{1}(x)D_{x_{n}}[\chi,op(\Lambda^{-\frac{1}{2}})]w \nonumber\\
&+& \mu T_{1}[\chi,op(\Lambda^{-\frac{1}{2}})]w
+\mu^{2}C_{0}[\chi,op(\Lambda^{-\frac{1}{2}})]w \nonumber\\
&=&op(\Lambda^{-\frac{1}{2}})Pw
+[P,op(\Lambda^{-\frac{1}{2}})]w+a_{1}+a_{2}+a_{3}+a_{4}+a_{5}.
\end{eqnarray}
Let us estimate $a_{1}$, $a_{2}$, $a_{3}$, $a_{4}$ and $a_{5}$.
Recalling that $[\chi,op(\Lambda^{-\frac{1}{2}})]\in
\mathcal{T}\mathcal {O}^{-\frac{3}{2}} $ and $\chi w=w$. Using the
fact that $[D_{x_{n}}, T_{k}]\in \mathcal{T}\mathcal{O}^{k}$ for all
$T_{k}\in \mathcal{T}\mathcal{O}^{k}$, we show
\begin{equation}\label{ttt2}
\norm{a_{1}}^{2}\leq
C\left(\norm{D_{x_{n}}^{2}op(\Lambda^{-\frac{3}{2}})w}^{2}+\norm{D_{x_{n}}op(\Lambda^{-\frac{3}{2}})w}^{2}+
\norm{op(\Lambda^{-\frac{3}{2}})w}^{2} \right)
\end{equation}
and
\begin{equation}\label{ttt3}
\norm{a_{3}}^{2}\leq C\left( \mu^{2}
\norm{D_{x_{n}}op(\Lambda^{-\frac{3}{2}})w}^{2}+\mu^{2}
\norm{op(\Lambda^{-\frac{3}{2}})w}^{2}\right).
\end{equation}
We have $R[\chi,op(\Lambda^{-\frac{1}{2}})]\in
\mathcal{T}\mathcal{O}^{\frac{1}{2}}$,
$T_{1}[\chi,op(\Lambda^{-\frac{1}{2}})]\in
\mathcal{T}\mathcal{O}^{-\frac{1}{2}}$ and \\ $
C_{0}[\chi,op(\Lambda^{-\frac{1}{2}})]\in
\mathcal{T}\mathcal{O}^{-\frac{3}{2}}$. Then we obtain
\begin{equation}\label{ttt4}
\norm{a_{2}}^{2}+\norm{a_{4}}^{2}+\norm{a_{5}}^{2}\leq C
\norm{op(\Lambda^{\frac{1}{2}})w}^{2}.
\end{equation}
Using the same computations made in the proof of Lemma \ref{ll} (cf
$t_{1}$, $t_{2}$ and $t_{3}$ of (\ref{ll2})), we show
\begin{equation}\label{ttt5}
\norm{[P,op(\Lambda^{-\frac{1}{2}})]w}^{2}\leq C\left( \norm{
op(\Lambda^{\frac{1}{2}})w}^{2}+\mu^{-1}\norm{D_{x_{n}}w}^{2}\right).
\end{equation}
Following (\ref{ttt1}), (\ref{ttt2}), (\ref{ttt3}), (\ref{ttt4}) and
(\ref{ttt5}), we obtain
\begin{equation}\label{ttt6}
\norm{Pu}^{2}\leq
C\left(\mu^{-1}\norm{Pw}^{2}+\norm{op(\Lambda^{\frac{1}{2}})w}^{2}+
\mu^{-1}\norm{D_{x_{n}}w}^{2}+
\mu^{-1}\norm{D_{x_{n}}^{2}op(\Lambda^{-1})w }^{2}\right).
\end{equation}
We have
$$op(b_{1})u=op(b_{1})\chi op(\Lambda^{-\frac{1}{2}})w=op(\Lambda^{-\frac{1}{2}})
op(b_{1})w+op(b_{1})[\chi, op(\Lambda^{-\frac{1}{2}})]w.$$ Recalling that $op(b_{1})\in
\mathcal{T}\mathcal{O}^{1}$, we obtain
\begin{equation}\label{ttt7}
\mu^{-1}\abs{op(b_{1})u}_{1}^{2}=\mu^{-1}\abs{op(\Lambda)op(b_{1})u}^{2}\leq
C\left( \mu^{-1}\abs {op(\Lambda^{\frac{1}{2}})op(b_{1})w
}^{2}+\mu^{-1}\abs {op(\Lambda^{\frac{1}{2}})w }^{2} \right).
\end{equation}
We have
$$op(b_{2})u=op(b_{2})\chi op(\Lambda^{-\frac{1}{2}})w=op(\Lambda^{-\frac{1}{2}})
op(b_{2})w+op(b_{2})[\chi, op(\Lambda^{-\frac{1}{2}})]w +[op(b_{2}),
op(\Lambda^{-\frac{1}{2}})]w.$$ Recalling that $op(b_{2})\in
D_{x_{n}} +\mathcal{T}\mathcal{O}^{1} $, we obtain
\begin{equation}\label{ttt8}
\mu\abs{op(b_{2})u}^{2}\leq C \left( \mu
\abs{op(\Lambda^{-\frac{1}{2}})op(b_{2})w}^{2}+\mu
\abs{op(\Lambda^{-\frac{1}{2}})w}^{2}+ \mu
\abs{D_{x_{n}}op(\Lambda^{-\frac{3}{2}})w}^{2} \right).
\end{equation}
Moreover, we have
$$\mu \abs{u}_{1,0,\mu}^{2}=\mu \abs{u}^{2}_{1}+ \mu\abs{D_{x_{n}}u}^{2}=\mu\abs{op(\Lambda)u}^{2}+\mu\abs{D_{x_{n}}u}^{2}. $$
We can write
$$op(\Lambda)u=op(\Lambda)\chi op(\Lambda^{-\frac{1}{2}})w= op(\Lambda^{\frac{1}{2}})w
 +op(\Lambda)[\chi, op(\Lambda^{-\frac{1}{2}})]w.$$
 Then $$ \mu \abs{op(\Lambda)u}^{2}\geq \mu \abs{ op(\Lambda^{\frac{1}{2}})w}^{2}-C \mu \abs{ op(\Lambda^{-\frac{1}{2}})w}^{2}\geq \mu \abs{ op(\Lambda^{\frac{1}{2}})w}^{2}- C\mu^{-1}\abs{ op(\Lambda^{\frac{1}{2}})w}^{2}.$$
 For $\mu$ large enough, we obtain
\begin{equation}\label{ttt9}
\mu \abs{op(\Lambda)u}^{2}\geq C\mu \abs{
op(\Lambda^{\frac{1}{2}})w}^{2}.
\end{equation}
By the same way, we prove, for $\mu$ large enough
\begin{equation}\label{ttt10}
\mu\abs{D_{x_{n}}u}^{2}\geq C\mu
\abs{D_{x_{n}}op(\Lambda^{-\frac{1}{2}})w}^{2}.
\end{equation}
Combining (\ref{ttt9}) and (\ref{ttt10}), we obtain
\begin{equation}\label{ttt11}
\mu\abs{u}_{1,0,\mu}^{2}\geq C\left( \mu \abs{
op(\Lambda^{\frac{1}{2}})w}^{2}+\mu
\abs{D_{x_{n}}op(\Lambda^{-\frac{1}{2}})w}^{2} \right).
\end{equation}
By the same way, we prove
\begin{equation}\label{ttt12}
\norm{op(\Lambda^{\frac{3}{2}})u}^{2}\geq \norm{op(\Lambda)w}^{2}-
C\norm{w}^{2},
\end{equation}
\begin{equation}\label{ttt13}
\norm{D_{x_{n}}op(\Lambda^{\frac{1}{2}})u}^{2}\geq
\norm{D_{x_{n}}w}^{2}-C\norm{op(\Lambda^{-1})D_{x_{n}}w}^{2}-C\norm{op(\Lambda^{-1})w}^{2}
\end{equation}
and
\begin{eqnarray}\label{ttt14}
\norm{
D_{x_{n}}^{2}op(\Lambda^{-\frac{1}{2}})u}^{2}\geq~~~~~~~~~~~~~~~~~~~~~~~~~
~~~~~~~~~~~~~~~~~~~~~~~~~~~~~~~~~~~~~~~~~~~~~~~~~~~~~\nonumber\\\\
\norm{D_{x_{n}}^{2}op(\Lambda^{-1})w}^{2}-C\norm{D_{x_{n}}^{2}op(\Lambda^{-2})w}^{2}-C
\norm{D_{x_{n}}op(\Lambda^{-2})w}^{2}-C\norm{op(\Lambda^{-2})w}^{2}.\nonumber
\end{eqnarray}
Combining (\ref{ttt12}), (\ref{ttt13}) and (\ref{ttt14}), we obtain
for $\mu$ large enough
\begin{eqnarray}\label{ttt15}
\norm{ D_{x_{n}}^{2}op(\Lambda^{-\frac{1}{2}})u}^{2}+
\norm{D_{x_{n}}op(\Lambda^{\frac{1}{2}})u}^{2}+
\norm{op(\Lambda^{\frac{3}{2}})u}^{2}~~~~~~~~~~~~~~~~~~~~~~~~~~~~~~~~~~~\nonumber\\\geq
C\left(\norm{D_{x_{n}}^{2}op(\Lambda^{-1})w}^{2}+\norm{D_{x_{n}}w}^{2}+\norm{op(\Lambda)w}^{2}
\right).
\end{eqnarray}
Combining (\ref{lll}), (\ref{ttt6}), (\ref{ttt7}), (\ref{ttt8}),
(\ref{ttt11}) and (\ref{ttt15}), we obtain (\ref{tt}), for $\mu$
large enough.\\
$~~~~~~~~~~~~~~~~~~~~~~~~~~~~~~~~~~~~~~~~~~~~~~~~~~~~~~~~~~~~~~
~~~~~~~~~~~~~~~~~~~~~~~~~~~~~~~~~~~~~~~~~~~~~~~~~~~~~~~~~~~~\square$
\section{Proof of Theorem \ref{t2}}
\hspace{5mm}This section is devoted to the proof of Theorem
\ref{t2}. \subsection{Study of the eigenvalues} \hspace{5mm} The
proof is based on a cutting argument related to the nature of the
roots of the polynomial $p_{j}(x,\xi',\xi_{n},\mu)$, $j=1,2$, in
$\xi_{n}$. On $x_{n}=0$, we note
$$ q_{1}(x',\xi',\mu)=q_{1,1}(0,x',\xi',\mu)=q_{1,2}(0,x',\xi',\mu).$$
Let us introduce the following micro-local regions
$$
\mathcal{E}^{+}_{1/_{2}}=\left\{(x,\xi',\mu)\in K\times\R^{n},\quad
q_{2,{1/_{2}}}+\frac{q_{1}^{2}}{(\frac{\partial
\varphi_{1/_{2}}}{\partial x_{n}})^{2}}>0 \right\},$$
$$ \mathcal{Z}_{1/_{2}}=\left\{(x,\xi',\mu)\in K\times\R^{n},\quad q_{2,{1/_{2}}}+\frac{q_{1}^{2}}{(\frac{\partial \varphi_{1/_{2}}}{\partial
x_{n}})^{2}}=0 \right\},$$
$$\mathcal{E}^{-}_{1/_{2}}=\left\{(x,\xi',\mu)\in K\times\R^{n},\quad q_{2,{1/_{2}}}+\frac{q_{1}^{2}}{(\frac{\partial \varphi_{1/_{2}}}{\partial
x_{n}})^{2}}<0 \right\}.$$ (Here and in the following the index
$1/_{2}$ using for telling $1$ or $2$).\\
 We decompose $p_{1/_{2}}(x,\xi,\mu)$ as a
polynomial in $\xi_{n}$. Then we have the following lemma describing
the various types of the roots of $p_{1/_{2}}$.
\begin{lemm}\label{l1}
We have the following
\begin{enumerate}
 \item  For $(x,\xi',\mu)\in \mathcal{E}^{+}_{1/_{2}}$, the roots of $p_{1/_{2}}$ denoted $z_{1/_{2}}^{\pm}$ satisfy $\pm\,\mbox{Im}\,z_{1/_{2}}^{\pm}>0$.
\item For $(x,\xi',\mu)\in \mathcal{Z}_{1/_{2}}$, one of the roots of $p_{1/_{2}}$ is real.
\item For $(x,\xi',\mu)\in \mathcal{E}^{-}_{1/_{2}}$, the roots of $p_{1/_{2}}$ are in the half- plane $\mbox{Im}\xi_{n}>0$
 if $\frac{\partial \varphi_{1/_{2}}}{\partial x_{n}}<0$ (resp. in the half-plane $\mbox{Im}\xi_{n}<0$ if $\frac{\partial \varphi_{1/_{2}}}{\partial x_{n}}>0$).
\end{enumerate}
 \end{lemm}
{\bf Proof.}\\
Using (\ref{e4}) and (\ref{e5}), we can write
\begin{equation}
\left \{\begin{array}{l}\label{e7}
p_{1}(x',\xi,\mu)=\left(\xi_{n}+i\mu \dfrac{\partial
\varphi_{1}}{\partial x_{n}}-i\alpha_{1}\right)\left(
\xi_{n}+i\mu \dfrac{\partial \varphi_{1}}{\partial x_{n}}+i\alpha_{1}\right),\\\\
p_{2}(x',\xi,\mu)=\left(\xi_{n}+i\mu \dfrac{\partial
\varphi_{2}}{\partial x_{n}}-i\alpha_{2}\right)\left( \xi_{n}+i\mu
\dfrac{\partial \varphi_{2}}{\partial x_{n}}+i\alpha_{2}\right),
 \end{array}
\right.
\end{equation}
where $\alpha_{j}\in \C$, $j=1,2$, defined by
\begin{equation}
\left \{\begin{array}{l}\label{e8}
\alpha_{1}^{2}(x',\xi',\mu)=\left(\mu \dfrac{\partial \varphi_{1}}{\partial x_{n}} \right)^{2}+q_{2,1}+2i\mu q_{1},\\\\
\alpha_{2}^{2}(x',\xi',\mu)=\left(\mu \dfrac{\partial
\varphi_{2}}{\partial x_{n}} \right)^{2}-\mu^{2}+q_{2,1}+2i\mu
q_{1}.
 \end{array}
\right.
\end{equation}
We set
\begin{equation}\label{e9}
z_{1/_{2}}^{\pm}=-i\mu \dfrac{\partial \varphi_{1/_{2}}}{\partial
x_{n}}\pm i\alpha_{1/_{2}},
\end{equation}
the roots of $p_{1/_{2}}$. The imaginary parts of the roots of
$p_{1/_{2}}$ are $$ -\mu \frac{\partial \varphi_{1/_{2}}}{\partial
x_{n}}- \mbox{Re}\, \alpha_{1/_{2}},\quad -\mu \frac{\partial
\varphi_{1/_{2}}}{\partial x_{n}}+ \mbox{Re}\, \alpha_{1/_{2}}.$$
The signs of the imaginary parts are opposite if $ \abs{\partial
\varphi_{1/_{2}}/\partial x_{n}}<\abs{\mbox{Re}\,\alpha_{1/_{2}}}$,
equal to the sign of $ -\partial \varphi_{1/_{2}}/\partial x_{n} $
if $ \abs{\partial \varphi_{1/_{2}}/\partial
x_{n}}>\abs{\mbox{Re}\,\alpha_{1/_{2}}}$ and one of the imaginary
parts is null if $ \abs{\partial \varphi_{1/_{2}}/\partial
x_{n}}=\abs{\mbox{Re}\,\alpha_{1/_{2}}}$. However the lines $
\mbox{Re}\,z=\pm \mu\, \partial \varphi_{1/_{2}}/\partial x_{n}$
 change by the application $ z\mapsto z'=z^{2}$ into the parabolic curve $ \mbox{Re}\,
  z'= \abs{\mu \,\partial \varphi_{1/_{2}}/\partial x_{n}}^{2}- \abs{\mbox{Im}\,z'}^{2}/4(\mu\, \partial \varphi_{1/_{2}}/\partial x_{n})^{2}$. Thus we obtain the
  lemma by replacing $z'$ by $\alpha_{1/_{2}}^{2}$. \\
$~~~~~~~~~~~~~~~~~~~~~~~~~~~~~~~~~~~~~~~~~~~~~~~~~~~~~~~~~~~~~~~~~~~~~~~~~~~~~~~~~~~~~~~~~~~~~~~~~~~~~~~~\square$
\begin{lemm}\label{l2}
If we assume that the function $\varphi$ satisfies the following
condition
\begin{equation}\label{e10}
 \left(\frac{\partial \varphi_{1}}{\partial x_{n}} \right)^{2}-  \left(\frac{\partial \varphi_{2}}{\partial x_{n}} \right)^{2}>1,
\end{equation}
then the following estimate holds
\begin{equation}\label{e11}
 q_{2,2}-\mu^{2}+\frac{q_{1}^{2}}{\left(\partial \varphi_{2}/\partial x_{n} \right)^{2} }>
 q_{2,1}+\frac{q_{1}^{2}}{\left(\partial \varphi_{1}/\partial x_{n} \right)^{2} }.
\end{equation}
\end{lemm}
{\bf Proof.} \\ Following (\ref{e5}), on $\{x_{n}=0\}$, we have
\begin{equation}\label{e12}
q_{2,2}(x,\xi',\mu)-q_{2,1}(x,\xi',\mu)=\left(\mu \frac{\partial
\varphi_{1}}{\partial x_{n}} \right)^{2}-\left(\mu\frac{\partial
\varphi_{2}}{\partial x_{n}} \right)^{2}.
\end{equation}
Using (\ref{e10}), we have (\ref{e11}).
$~~~~~~~~~~~~~~~~~~~~~~~~~~~~~~~~~~~~~~~~~~~~~~~~~~~~~~~~~~~~~~\square$
\begin{rema}\label{r1}
The result of this lemma imply that $\mathcal{E}^{+}_{1}\subset
\mathcal{E}^{+}_{2}$.
\end{rema}
\subsection{Estimate in $\mathcal{E}^{+}_{1}$} \hspace{5mm}In this
part we study the problem in the elliptic region $\mathcal{E}^{+}$.
In this region we can inverse the operator and use the Calderon
projectors. Let $\chi^{+}(x,\xi',\mu)\in \mathcal{T}S_{\mu}^{0}$
such that in the support of $\chi^{+}$ we have $ q_{2,1}+
\frac{q_{1}^{2}}{(\partial \varphi_{1}/\partial x_{n})^{2}}\geq
\delta > 0$. Then we have the following partial estimate.
\begin{prop}\label{p1}
There exists a constant $C>0$ and $\mu_{0}>0$ such that for any $\mu
\geq \mu_{0}$, we have
\begin{equation}\label{e17}
 \mu^{2}\norm{op(\chi^{+})u}_{1,\mu}^{2}\leq C\left( \norm{P(x,D,\mu)u}^{2}+\norm{u}_{1,\mu}^{2}+
 \mu\abs{u}_{1,0,\mu}^{2}\right),
\end{equation}
for any $u\in C_{0}^{\infty}(\overline{\Omega}_{2})$.\par If we suppose
moreover that $ \varphi$ satisfies (\ref{e10}) then the following
estimate holds
\begin{equation}\label{e18}
\mu \abs{op(\chi^{+})u}_{1,0,\mu}^{2}\leq C\left(
\norm{P(x,D,\mu)u}^{2}+\mu^{-1} \abs{op(b_{1})u}^{2}_{1}+\mu
\abs{op(b_{2})u}^{2}+
 \norm{u}_{1,\mu}^{2}+\mu^{-2}\abs{u}^{2}_{1,0,\mu}
 \right),
\end{equation}
for any $u\in C_{0}^{\infty}(\overline{\Omega}_{2})$ and $b_{j}$,
$j=1,2$, defined in (\ref{e16}).
\end{prop}
{\bf Proof}\\
Let $\tilde{u}=op(\chi^{+})u$. Then we get

\begin{equation}
\left \{\begin{array}{lc}\label{e19}
P\tilde{u}=\tilde{f}&\mbox{in}\, \left\{x_{n}>0\right\},\\
op (b_{1})\tilde{u}=\tilde{u}_{0}|_{x_{n}=0}-i\mu \tilde{v}_{0}|_{x_{n}=0}=\tilde{e}_{1}&\mbox{on}\, \left\{x_{n}=0\right\},\\
op(b_{2})\tilde{u}=\left(D_{x_{n}}+ i\mu \frac{\partial
\varphi_{1}}{\partial x_{n}}\right)\tilde{u}_{0}|_{x_{n}=0}+
\left(D_{x_{n}}+i\mu \frac{\partial \varphi_{2}}{\partial
x_{n}}\right)\tilde{v}_{0}|_{x_{n}=0}=\tilde{e}_{2}&\mbox{on}\,\left\{x_{n}=0\right\},
\end{array}
\right.
\end{equation}
with $ \tilde{f}=op(\chi^{+})f+ \left[ P, op(\chi^{+})\right]u $.
Since $ \left[ P, op(\chi^{+})\right]\in
(\mathcal{T}\mathcal{O}^{0})D_{x_{n}}+ \mathcal{T}\mathcal{O} ^{1}$,
we have
\begin{equation}\label{n0}
\| \tilde{f}\|_{L^{2}}^{2}\leq C \left(
\norm{P(x,D,\mu)u}_{L^{2}}^{2}+\norm{u}_{1,\mu}^{2}\right)
\end{equation}
and $\tilde{e}_{1}=op(\chi^{+})e_{1}$ satisfying
\begin{equation}\label{n'0}
\abs{\tilde{e}_{1}}_{1}^{2}\leq C\abs{e_{1}}_{1}^{2}
\end{equation}
 and \\\\
$\tilde{e}_{2}= \left[ (D_{x_{n}}+i\mu \frac{\partial
\varphi_{1}}{\partial x_{n}}), op(\chi^{+})\right]u_{0}|_{x_{n}=0}+
\left[ (D_{x_{n}}+i\mu \frac{\partial \varphi_{2}}{\partial x_{n}}),
op(\chi^{+})\right]v_{0}|_{x_{n}=0}+ op(\chi^{+})e_{2}.$ \\\\ Since
$ [D_{x_{n}}, op(\chi^{+})]\in \mathcal{T}\mathcal{O}^{0}$, we have
\begin{equation}\label{n00}
| \tilde{e}_{2}|^{2}\leq C \left( |u|^{2}+ |e_{2}|^{2}\right).
\end{equation}
 Let $ \underline{\tilde{u}}$ the extension of $\tilde{u}$ by $0$
in $x_{n}<0$. According to (\ref{e4}), (\ref{e5}) and (\ref{e6}), we
obtain, by noting $\partial \varphi/\partial
x_{n}=\mbox{diag}\left(\partial \varphi_{1}/\partial x_{n}, \partial
\varphi_{2}/\partial x_{n}\right) $, $
\gamma_{j}(\tilde{u})=\,^t\left( D_{x_{n}}^{j}(\tilde{u}_{0})\mid
_{x_{n}=0^{+}},D_{x_{n}}^{j}(\tilde{v}_{0})\mid
_{x_{n}=0^{+}}\right)$, $j=0,1$ and $\delta^{(j)}=\left(d/d x_{n}
\right)^{j}\left(\delta_{x_{n}=0}\right)$,
\begin{equation}\label{e20}
P\underline{\tilde{u}}=
\underline{\tilde{f}}-\gamma_{0}(\tilde{u})\otimes
\delta'+\frac{1}{i} \left(\gamma_{1}(\tilde{u})+2i\mu \frac{\partial
\varphi}{\partial x_{n}} \right)\otimes \delta
\end{equation}
Let $ \chi(x,\xi,\mu)\in S_{\mu}^{0}$ equal to $ 1$ for sufficiently
large $|\xi|+\mu$ and in a neighborhood of supp($ \chi^{+}$) and
satisfies that in the support of $\chi$ we have $p$ is elliptic.
These conditions are compatible  from the choice made for supp($
\chi^{+}$) and Remark \ref{r1}. Let $m$ large enough chosen later,
by the ellipticity of $p$ on supp$( \chi)$ there exists $E=op(e)$ a
parametrix of $P$. We recall that $e\in S_{\mu}^{-2}$, of the form
$e(x,\xi,\mu)=\sum_{j=0}^{m}e_{j}(x,\xi,\mu)$, where $ e_{0}=\chi
p^{-1}$ and $e_{j}=\mbox{diag}(e_{j,1},e_{j,2})\in S_{\mu}^{-2-j}$
such that $e_{j,1}$ and $e_{j,2}$ are rational fractions in
$\xi_{n}$. Then we have
\begin{equation}\label{e21}
EP=op(\chi)+ R_{m},\quad R_{m}\in \mathcal{O}^{-m-1}.
\end{equation} Following (\ref{e20}) and
(\ref{e21}), we obtain
\begin{equation}\left\{\begin{array}{l}\label{e22}
\underline{\tilde{u}}=E\underline{\tilde{f}}+E\left[  -h_{1}\otimes \delta'+\dfrac{1}{i}h_{0}\otimes \delta\right] +w_{1},\\\\
h_{0}=\gamma_{1}(\tilde{u})+2i\mu \dfrac{\partial \varphi}{\partial x_{n}}\gamma_{0}(\tilde{u}),\quad h_{1}=\gamma_{0}(\tilde{u}),\\\\
w_{1}=
\left(\mbox{Id}-op(\chi)\right)\underline{\tilde{u}}-R_{m}\underline{\tilde{u}}.
\end{array}
\right.
\end{equation}
Using the fact that supp$(1-\chi)\cap
\mbox{supp}(\chi^{+})=\varnothing$ and symbolic calculus (See Lemma
2.10 in \cite{LRR}), we have $
\left(\mbox{Id}-op(\chi)\right)op(\chi^{+})\in \mathcal{O}^{-m}$,
then we obtain
\begin{equation}\label{n01}
 \norm{w_{1}}_{2,\mu}^{2}\leq C \mu^{-2}\norm{u}_{L^{2}}^{2}.
\end{equation}
 Now, let us look at this term $ E\left[  -h_{1}\otimes \delta'+\dfrac{1}{i}h_{0}\otimes \delta\right] $.
  For $x_{n}>0$, we get
$$\left \{\begin{array}{l}
 E\left[  -h_{1}\otimes \delta'+\dfrac{1}{i}h_{0}\otimes
 \delta\right]=\hat{T}_{1}h_{1}+\hat{T}_{0}h_{0},\\\\
 \hat{T}_{j}(h)=\left(\dfrac{1}{2\pi}\right)^{n-1}\dint
 e^{i(x'-y')\xi'}\hat{t}_{j}(x,\xi',\mu)h(y')dy'd\xi'=op
 (\hat{t}_{j})h\\\\
 \hat{t}_{j}=\dfrac{1}{2\pi
 i}\dint_{\gamma}e^{ix_{n}\xi_{n}}e(x,\xi,\mu)\xi_{n}^{j}d\xi_{n}
\end{array}
\right.
 $$
where $\gamma$ is the union of the segment $\{\xi_{n}\in \R,\,
|\xi_{n}|\leq c_{0}\sqrt{|\xi'|^{2}+\mu^{2}} \}$ and the half circle
$\{\xi_{n}\in \C,\, |\xi_{n}|= c_{0}\sqrt{|\xi'|^{2}+\mu^{2}},\,Im
\xi_{n}>0 \} $, where the constant $c_{0}$ is chosen sufficiently
large so as to have the roots $z_{1}^{+}$ and $z_{2}^{+}$ inside the
domain with boundary $\gamma$ (If $c_{0}$ is large enough, the
change of contour $\R \longrightarrow \gamma$ is possible because
the symbol $e(x,\xi,\mu)$ is holomorphic for large $ |\xi_{n}|$;
$\xi_{n}\in C$). In particular we have in $x_{n}\geq 0$
\begin{equation}\label{n1}
\abs{\partial_{x_{n}}^{k}\partial_{x'}^{\alpha}\partial_{\xi'}^{\beta}\hat{t}_{j}}\leq
C_{\alpha,\beta,k} \langle \xi',\mu\rangle^{j-1-\abs{\beta}+k},\quad
j=0,1.
\end{equation}
We now choose $\chi_{1}(x,\xi',\mu)\in \mathcal{T}S_{\mu}^{0}$,
satisfying the same requirement as $\chi^{+}$, equal to $1$ in a
neighborhood of supp$(\chi^{+})$ and such that the symbol $\chi$ be
equal to $1$ in a neighborhood of supp$(\chi_{1})$. We set
$t_{j}=\chi_{1}\hat{t}_{j}$, $j=0,1$. Then we obtain
\begin{equation}\label{n2}
\underline{\tilde{u}}=E\underline{\tilde{f}}+op(t_{0})h_{0}+
op(t_{1})h_{1}+w_{1}+w_{2}
\end{equation}
where $ w_{2}=
op((1-\chi_{1})\hat{t}_{0})h_{0}+op((1-\chi_{1})\hat{t}_{1})h_{1}$.
By using the composition formula of tangential operator, estimate
(\ref{n1}), the fact that supp$(1-\chi_{1})\cap
\mbox{supp}(\chi^{+})=\varnothing$ and the following trace formula
\begin{equation}\label{n3}
|\gamma_{0}(u)|_{j}\leq C\mu^{-\frac{1}{2}}\|u\|_{j+1,\mu},\quad
j\in \N,
\end{equation}we obtain
\begin{equation}\label{n4}
\|w_{2}\|_{2,\mu}^{2}\leq C\mu^{-2}\left(
\|u\|_{1,\mu}^{2}+|u|_{1,0,\mu}^{2}\right).
\end{equation}
Since $\chi=1$ in the support of $\chi_{1}$, we have $e(x,\xi,\mu)$
is meromorphic w.r.t $\xi_{n}$ in the support of $\chi_{1}$.
$z^{+}_{1/_{2}}$ are in $\mbox{Im}\xi_{n}\geq
c_{1}\sqrt{|\xi'|^{2}+\mu^{2}} $ $(c_{1}>0)$. If $c_{1}$ is small
enough we can choose $ \gamma_{1/_{2}}$ in $\mbox{Im}\xi_{n}\geq
\frac{c_{1}}{2}\sqrt{|\xi'|^{2}+\mu^{2}} $ and we can write
\begin{equation}\label{n5}
t_{j}=\mbox{diag}(t_{j,1}, t_{j,2}),\quad
t_{j,1/_{2}}(x,\xi',\mu)=\chi_{1}(x,\xi',\mu)\frac{1}{2\pi
i}\int_{\gamma_{1/_{2}}}e^{ix_{n}\xi_{n}}e_{1/_{2}}(x,\xi,\mu)\xi_{n}^{j}d\xi_{n},\quad
j=0,1.
\end{equation}
Then there exists $ c_{2}>0$ such that in $x_{n}\geq 0$, we obtain
\begin{equation}\label{n7}
\abs{\partial_{x_{n}}^{k}\partial_{x'}^{\alpha}\partial_{\xi'}^{\beta}t_{j}}\leq
C_{\alpha,\beta,k}e^{-c_{2}x_{n}\langle \xi',\mu\rangle} \langle
\xi',\mu\rangle^{j-1-|\beta|+k}.
\end{equation}
In particular, we have
$e^{c_{2}x_{n}\mu}(\partial_{x_{n}}^{k})t_{j}$ is bounded in $
\mathcal{T}S_{\mu}^{j-1+k}$ uniformly w.r.t $ x_{n}\geq 0$. Then
$$ \norm{\partial_{x'}op(t_{j})h_{j}}_{L^{2}}^{2}+ \norm{op(t_{j})h_{j}}_{L^{2}}^{2}\leq
C \int_{x_{n}>0}e^{-2 c_{2}x_{n}\mu} \abs{ op(
e^{c_{2}x_{n}\mu}t_{j})h_{j}}_{1}^{2}(x_{n})dx_{n}\leq C \mu^{-1}
|h_{j}|_{j}^{2}$$ and
$$\norm{\partial_{x_{n}}op(t_{j})h_{j}}_{L^{2}}^{2}\leq C \int_{x_{n}>0}e^{-2 c_{2}x_{n}\mu}
\abs{ op( e^{c_{2}x_{n}\mu}
\partial_{x_{n}}t_{j})h_{j}}_{L^{2}}^{2}(x_{n})dx_{n}\leq C \mu^{-1}
|h_{j}|_{j}^{2}. $$
 Using the fact that $
h_{0}=\gamma_{1}(\tilde{u})+2i\mu \frac{\partial \varphi}{\partial
x_{n}}\gamma_{0}(\tilde{u})$ and $ h_{1}=\gamma_{0}(\tilde{u})$, we
obtain
\begin{equation}\label{n7}
\norm{op(t_{j})h_{j}}_{1,\mu}^{2}\leq C \mu^{-1}|u|_{1,0,\mu}^{2}.
\end{equation}
From (\ref{n2}) and estimates (\ref{n0}), (\ref{n01}), (\ref{n4})
and (\ref{n7}), we obtain (\ref{e17}).\\
It remains to proof (\ref{e18}).  We recall that, in
supp($\chi_{1}$), we have
$$ e_{0}=\mbox{diag}\left(e_{0,1},e_{0,2}\right)=
\mbox{diag}\left(\frac{1}{p_{1}},\frac{1}{p_{2}}\right)=
\mbox{diag}\left( \frac{1}{(\xi_{n}-z_{1}^{+})(\xi_{n}-z_{1}^{-})},
 \frac{1}{(\xi_{n}-z_{2}^{+})(\xi_{n}-z_{2}^{-})}\right).
$$
Using the residue formula, we obtain
\begin{equation}\label{n6}
e^{-ix_{n}z_{1/_{2}}^{+}}t_{j,1/_{2}}=\chi_{1}\frac{
(z_{1/_{2}}^{+})^{j}}{ z_{1/_{2}}^{+}-z_{1/_{2}}^{-}}+
\lambda_{1/_{2}},\quad j=0,1,\quad \lambda_{1/_{2}}\in
\mathcal{T}S_{\mu}^{-2+j}.
\end{equation}
Taking the traces of (\ref{n2}), we obtain
\begin{equation}\label{e28}
\gamma_{0}(\tilde{u})=op(c)\gamma_{0}(\tilde{u})+
op(d)\gamma_{1}(\tilde{u})+w_{0},
\end{equation}
where $ w_{0}=\gamma_{0}(E \underline{\tilde{f}}+w_{1}+w_{2})$
satisfies, according to the trace formula (\ref{n3}), the estimates
(\ref{n0}), (\ref{n01}) and (\ref{n4}), the following estimate
\begin{equation}\label{n9}
\mu\abs{w_{0}}_{1}^{2}\leq C\left(
\norm{P(x,D,\mu)u}^{2}+\norm{u}_{1,\mu}^{2}+\mu^{-2}\abs{u}_{1,0,\mu}^{2}\right)
\end{equation}
and following (\ref{n7}), $c$ and $d$ are two tangential symbols of
order respectively $0$ and $-1$ given by
$$\begin{array}{lll}
c_{0}=\mbox{diag}(c_{0,1},c_{0,2})&\mbox{with}&c_{0,1/_{2}}=
-\left(\chi_{1}\dfrac{z_{1/_{2}}^{-}}
{z_{1/_{2}}^{+}- z_{1/_{2}}^{-}}\right),\\\\
d_{-1}=\mbox{diag}(d_{-1,1},d_{-1,2})&\mbox{with}&d_{-1,1/_{2}}=
\left(\chi_{1}\dfrac{1} {z_{1/_{2}}^{+}- z_{1/_{2}}^{-}}\right).
\end{array}$$
Following (\ref{e19}), the transmission conditions give
\begin{equation}\left\{\begin{array}{l}\label{e32}
\gamma_{0}(\tilde{u}_{0})-i\mu \gamma_{0}(\tilde{v}_{0})=\tilde{e}_{1}\\\\
\gamma_{1}(\tilde{u}_{0})+ \gamma_{1}(\tilde{v}_{0})+i\mu
\frac{\partial \varphi_{1}}{\partial
x_{n}}\gamma_{0}(\tilde{u}_{0})+i\mu \frac{\partial
\varphi_{2}}{\partial x_{n}}\gamma_{0}(\tilde{v}_{0})=\tilde{e}_{2}.
\end{array}
\right.
\end{equation}
We recall that $ \tilde{u}= (\tilde{u}_{0},\tilde{v}_{0})$,
combining (\ref{e28}) and (\ref{e32}) we show that
\begin{equation}\label{e33}
op(k)\,^t\left( \gamma_{0}(\tilde{u}_{0}),
\gamma_{0}(\tilde{v}_{0}),
 \Lambda^{-1}\gamma_{1}(\tilde{u}_{0}), \Lambda^{-1}\gamma_{1}(\tilde{v}_{0})\right)=
 w_{0}+\frac{1}{\mu}op \left(
                         \begin{array}{c}
                           0 \\
                           0 \\
                           1 \\
                           0 \\
                         \end{array}
                       \right)\tilde{e}_{1}+op \left(
                         \begin{array}{c}
                           0 \\
                           0 \\
                           0 \\
                           1 \\
                         \end{array}
                       \right)\Lambda^{-1}\tilde{e}_{2}
 ,
\end{equation}
where $k$ is a $ 4\times 4$ matrix, with principal symbol defined by
$$k_{0}+\frac{1}{\mu}r_{0}=\left( \begin{array}{cccc}
1-c_{0,1}&0&-\Lambda\, d_{-1,1}&0\\\\
0&1-c_{0,2}&0&-\Lambda \, d_{-1,2}\\\\
0&-i&0&0\\\\
i\mu\Lambda^{-1}\dfrac{\partial \varphi_{1}}{\partial
x_{n}}&i\mu\Lambda^{-1}\dfrac{\partial \varphi_{2}}{\partial
x_{n}}&1&1
                                               \end{array}
 \right)+\frac{1}{\mu}r_{0},
  $$
  where $r_{0}$ is a tangential symbol of order $0$.\\
  We now choose $\chi_{2}(x,\xi',\mu)\in \mathcal{T}S_{\mu}^{0}$,
  satisfying the same requirement as $\chi^{+}$, equal to $1$ in a
  neighborhood of supp$(\chi^{+})$ and such that the symbol $
  \chi_{1}$ be equal to $1$ in a neighborhood of supp$(\chi_{2})$.
  In supp$(\chi_{2})$, we obtain
  $$ k_{0}|_{ \mbox{supp}(\chi_{2})}=\left( \begin{array}{cccc}
\dfrac{z_{1}^{+}}{z_{1}^{+}-z_{1}^{-}}&0&-\dfrac{\Lambda}{z_{1}^{+}-z_{1}^{-}}&0\\\\
0&\dfrac{z_{2}^{+}}{z_{2}^{+}-z_{2}^{-}}&0&-\dfrac{\Lambda}{z_{2}^{+}-z_{2}^{-}}\\\\
0&-i&0&0\\\\
i\mu\Lambda^{-1}\dfrac{\partial \varphi_{1}}{\partial
x_{n}}&i\mu\Lambda^{-1}\dfrac{\partial \varphi_{2}}{\partial x_{n}}
&1&1
                                               \end{array}
 \right).$$
Then, following (\ref{e9}),
$$ \mbox{det}(k_{0})|_{ \mbox{supp}(\chi_{2})}=-\left( z_{1}^{+}-z_{1}^{-}\right)^{-1}\left( z_{2}^{+}-z_{2}^{-}\right)^{-1}\Lambda\,
\alpha_{1}.$$ To prove that there exists $c>0$ such that $\abs{
\mbox{det}(k_{0})|_{ \mbox{supp}(\chi_{2})}}\geq c $, by homogeneity
it suffices to prove that $ \mbox{det}(k_{0})|_{
\mbox{supp}(\chi_{2})}\neq 0$ if $ \abs{\xi'}^{2}+\mu^{2}=1$.\\
 If we suppose that $\mbox{det}(k_{0})|_{
\mbox{supp}(\chi_{2})}=0 $, we obtain $\alpha_{1}=0$ and then
$\alpha_{1}^{2}=0$.\\ Following (\ref{e8}),we obtain
$$
q_{1}=0 \quad\mbox{and}\quad \left(\mu \frac{\partial
\varphi_{1}}{\partial x_{n}}\right)^{2}+q_{2,1}=0.
$$ Combining with the fact that
$q_{2,1}+\frac{q_{1}^{2}}{\left( \partial \varphi_{1}/\partial x_{n}
\right)^{2} }>0$, we obtain
$$-\left( \mu \frac{\partial \varphi_{1}}{\partial
x_{n}}\right)^{2}>0.$$
 Therefore $ \mbox{det}(k_{0})|_{ \mbox{supp}(\chi_{2})}\neq 0$. It follows that, for large $\mu$, $k=k_{0}+\frac{1}{\mu}r_{0}$ is elliptic
in supp($\chi_{2}$). Then there exists $l\in
\mathcal{T}S_{\mu}^{0}$, such that
$$op(l)op(k)=op(\chi_{2})+ \tilde{R}_{m},$$
with $ \tilde{R}_{m} \in \mathcal {T}\mathcal{O}^{-m-1}$, for $m$
large. This yields
$$\begin{array}{l} ^t\left( \gamma_{0}(\tilde{u}_{0}),
\gamma_{0}(\tilde{v}_{0}),
 \Lambda^{-1}\gamma_{1}(\tilde{u}_{0}),
 \Lambda^{-1}\gamma_{1}(\tilde{v}_{0})\right)=
op(l)w_{0}+\frac{1}{\mu}op(l)op \left(
                         \begin{array}{c}
                           0 \\
                           0 \\
                           1 \\
                           0 \\
                         \end{array}
                       \right)\tilde{e}_{1}+ op(l)op \left(
                         \begin{array}{c}
                           0 \\
                           0 \\
                           0 \\
                           1 \\
                         \end{array}
                       \right)\Lambda^{-1}\tilde{e}_{2}\\
                       ~~~~~~~~~~~~~~~~~~~~~~~~~~~~~~~~~~~~~~~~+(op(1-\chi_{2})-\tilde{R}_{m})
^t\left( \gamma_{0}(\tilde{u}_{0}), \gamma_{0}(\tilde{v}_{0}),
 \Lambda^{-1}\gamma_{1}(\tilde{u}_{0}), \Lambda^{-1}\gamma_{1}(\tilde{v}_{0})\right)
 .\end{array}$$
 Since supp$(1-\chi_{2})\cap \mbox{supp}(\chi^{+})=\varnothing$ and
 by using (\ref{n9}), we obtain
 $$ \mu|\tilde{u}|_{1,0,\mu}^{2}\leq C\left( \mu^{-1} |\tilde{e}_{1}|_{1}^{2}+\mu
 |\tilde{e}_{2}|+ \norm{P(x,D,\mu)u}_{L^{2}}^{2}+\norm{u}_{1,\mu}^{2}+\mu^{-2}
 \abs{u}_{1,0,\mu}^{2}\right).$$
 From estimates (\ref{n'0}) and (\ref{n00}) and the trace formula (\ref{n3}), we
 obtain (\ref{e18}).\\$~~~~~~~~~~~~~~~~~~~~~~~~~~~~~~~~~~~~~~~~~~~~~~~~~~~~~~~~~~~~~~~~~~~~~~~~~~~~~~~~~~~~~~
 ~~~~~~~~~~~~~~~~~~~~~~~~~~~~~~~~~~~~~~~~~~~~~~~~\square$
\subsection{Estimate in $\mathcal{Z}_{1}$} \hspace{5mm}The aim of
this part is to prove the estimate in the region $\mathcal {Z}_{1}$.
In this region, if $ \varphi$ satisfies (\ref{e10}), the symbol $
p_{1}(x,\xi,\mu)$ admits a real roots and $p_{2}(x,\xi,\mu)$ admits
two roots $z_{2}^{\pm}$ satisfy $\pm\mbox{ Im} (z_{2}^{\pm})>0$. Let
$ \chi^{0}(x,\xi',\mu)\in \mathcal{T}\mathcal{S}_{\mu}^{0}$ equal to
$ 1$ in $\mathcal{Z}_{1}$ and such that in the support of $\chi^{0}$
we have $q_{2,2}-\mu^{2}+ \frac{q_{1}^{2}}{ (\partial
\varphi_{2}/\partial x_{n})^{2}}\geq \delta>0$. Then we have the
following partial estimate.
\begin{prop}\label{p2} There exists constants $C>0$ and
$\mu_{0}>0$ such that for any $\mu\geq\mu_{0}$ we have the following
estimate
\begin{equation}\label{z1}
\mu\norm{op(\chi^{0})u}^{2}_{1,\mu}\leq
C\left(\norm{P(x,D,\mu)u}^{2}+\mu\abs{u}^{2}_{1,0,\mu}+\norm{u}^{2}_{1,\mu}\right),
\end{equation}
for any $u\in C_{0}^{\infty}(\overline{\Omega}_{2})$.\par If we assume
moreover that $ \varphi$ satisfies (\ref{e10}) then we have
\begin{equation}\label{z2}
\mu\abs{op(\chi^{0})u}^{2}_{1,0,\mu}\leq
C\left(\norm{P(x,D,\mu)u}^{2}+\mu^{-1} \abs{op(b_{1})u}^{2}_{1}+\mu
\abs{op(b_{2})u}^{2}+
 \norm{u}_{1,\mu}^{2}+\mu^{-2}\abs{u}^{2}_{1,0,\mu}
 \right),
\end{equation}
for any $u\in C_{0}^{\infty}(\overline{\Omega}_{2})$ and $b_{j}$,
$j=1,2$, defined in (\ref{e16}).
\end{prop}
 \subsubsection {Preliminaries}\hspace{5mm}
Let $u\in C_{0}^{\infty}(K)$, $ \tilde{u}=op(\chi^{0})u$ and $P$ the
differential operator with principal symbol given by
$$ p(x,\xi,\mu)=\xi_{n}^{2}+q_{1}(x,\xi',\mu)\xi_{n}+q_{2}(x,\xi',\mu)$$
where $q_{j}=\mbox{diag}(q_{j,1},q_{j,2})$, $j=1,2$. Then we have
the following system
\begin{equation}\left\{\begin{array}{ll}\label{z3}
P\tilde{u}=\tilde{f}&\mbox{in}\,\{ x_{n}>0\},\\\\
B\tilde{u}=\tilde{e}=(\tilde{e}_{1},\tilde{e}_{2})&\mbox{on}\,\{x_{n}=0\},
\end{array}
\right.\end{equation} where
$\tilde{f}=op(\chi^{0})f+\left[P,op(\chi^{0})\right]u$. Since $
\left[ P, op(\chi^{0})\right]\in
(\mathcal{T}\mathcal{O}^{0})D_{x_{n}}+ \mathcal{T}\mathcal{O} ^{1}$,
we have
\begin{equation}\label{f1}
\| \tilde{f}\|_{L^{2}}^{2}\leq C \left(
\norm{P(x,D,\mu)u}_{L^{2}}^{2}+\norm{u}_{1,\mu}^{2}\right),
\end{equation}
$B$ defined in (\ref{e16}) and  $\tilde{e}_{1}=op(\chi^{0})e_{1}$
satisfying
\begin{equation}\label{f2}
\abs{\tilde{e}_{1}}_{1}^{2}\leq C\abs{e_{1}}_{1}^{2}
\end{equation}
and \\\\
$\tilde{e}_{2}= \left[ (D_{x_{n}}+i\mu \frac{\partial
\varphi_{1}}{\partial x_{n}}), op(\chi^{0})\right]u_{0}|_{x_{n}=0}+
\left[ (D_{x_{n}}+i\mu \frac{\partial \varphi_{2}}{\partial x_{n}}),
op(\chi^{0})\right]v_{0}|_{x_{n}=0}+ op(\chi^{0})e_{2}.$ \\\\
Since $ [D_{x_{n}}, op(\chi^{+})]\in \mathcal{T}\mathcal{O}^{0}$, we
have
\begin{equation}\label{f3}
| \tilde{e}_{2}|^{2}\leq C \left( |u|^{2}+ |e_{2}|^{2}\right).
\end{equation}
\par Let us reduce the problem (\ref{z3}) to a first order system.
Put $v=^t\left(\langle D',\mu\rangle
\tilde{u},D_{x_{n}}\tilde{u}\right)$. Then we obtain the following
system
\begin{equation}\left\{\begin{array}{ll}\label{z4}
D_{x_{n}}v-op( \mathcal{P})v=F&\mbox{in}\,\{x_{n}>0\},\\\\
op(\mathcal{B})v=(\frac{1}{\mu}\Lambda\tilde{e}_{1},\tilde{e}_{2})&\mbox{on}\,\{x_{n}=0\},
\end{array}
\right.
\end{equation}
where $\mathcal{P}$ is a $4\times4$ matrix, with principal symbol
defined by
$$\mathcal{P}_{0}= \left( \begin{array}{cc}
0&\Lambda\,\mbox{Id}_{2}\\
\Lambda^{-1}q_{2}&-q_{1}\\
                                               \end{array}
 \right),\quad \left( \Lambda= \langle \xi',\mu\rangle=\left( \abs{\xi'}^{2}+\mu^{2}\right)^{\frac{1}{2}}\right),$$
  $ \mathcal{B}$ is a tangential symbol of order $0$, with principal symbol given by
  $$ \mathcal{B}_{0}+\frac{1}{\mu}r_{0}=\left( \begin{array}{cccc}
0&-i&0&0\\
i\mu\Lambda^{-1}\frac{\partial \varphi_{1}}{\partial
x_{n}}&i\mu\Lambda^{-1}\frac{\partial \varphi_{2}}{\partial
x_{n}}&1&1
\end{array}
 \right)+\frac{1}{\mu}r_{0} $$
 ($r_{0}$ a tangential symbol of order $0$) and $F=^t(0,\tilde{f})$.
 \par For a fixed $(x_{0},\xi'_{0},\mu_{0})$ in supp$\chi_{0}$, the
 generalized eigenvalues of the matrix $\mathcal{P}$ are the zeroes
 in $\xi_{n}$ of $p_{1}$ and $p_{2}$ i.e $z_{1}^{\pm}=-i\mu\frac{\partial \varphi_{1}}{\partial
x_{n}}\pm i\alpha_{1} $ and $z_{2}^{\pm}=-i\mu\frac{\partial
\varphi_{2}}{\partial x_{n}}\pm i\alpha_{2}$ with $\pm
\mbox{Im}(z_{2}^{\pm})>0$ and $ z_{1}^{+}\in {\R}$.
\par We note $s(x,\xi',\mu)=( s_{1}^{-},s_{2}^{-},s_{1}^{+},s_{2}^{+})$ a basis of the
generalized eigenspace of $\mathcal{P}(x_{0},\xi'_{0},\mu_{0})$
corresponding to eigenvalues with positive or negative imaginary
parts. $s_{j}^{\pm}(x,\xi',\mu)$, $j=1,2$ is a $C^{\infty}$ function
on a conic neighborhood of $(x_{0},\xi'_{0},\mu_{0})$ of a degree
zero in $(\xi',\mu)$. We denote $op(s)(x,D_{x'},\mu)$ the
pseudo-differential operator associated to the principal symbol $
s(x,\xi',\mu)=\left(s_{1}^{-}(x,\xi',\mu),s_{2}^{-}(x,\xi',\mu),s_{1}^{+}(x,\xi',\mu),
s_{2}^{+}(x,\xi',\mu)\right)$.\\
Let $ \hat{\chi}(x,\xi',\mu)\in \mathcal{T}S_{\mu}^{0}$ equal to $1$
in a conic neighborhood of $(x_{0},\xi'_{0},\mu_{0})$ and in a
neighborhood of supp$ (\chi^{0})$ and satisfies that in the support
of $\hat{\chi}$, $ s$ is elliptic. Then there exists $ n\in
\mathcal{T}S_{\mu}^{0}$, such that
$$ op(s)op(n)=op( \hat{\chi})+ \hat{R}_{m},$$
with $ \hat{R}_{m}\in \mathcal{T}\mathcal{O}^{-m-1}$, for $m$
large.\\ Let $V=op(n)v$. Then we have the following system
\begin{equation}\left\{\begin{array}{ll}\label{f4}
D_{x_{n}}V= GV +AV+F_{1}&\mbox{in}\,\{x_{n}>0\},\\\\
op(\mathcal{B}_{1})V=(\frac{1}{\mu}\Lambda\tilde{e}_{1},\tilde{e}_{2})
+v_{1}&\mbox{on}\,\{x_{n}=0\},
\end{array}
\right.
\end{equation}
where $G=op(n)op(\mathcal{P})op(s)$, $A=
\left[D_{x_{n}},op(n)\right]op(s)$,\\
$F_{1}=op(n)F+op(n)op(\mathcal{P})(op(1- \hat{\chi})- \hat{R}_{m})v+
\left[D_{x_{n}},op(n)\right](op(1- \hat{\chi})- \hat{R}_{m})v$, $op(
\mathcal{B}_{1})= op(\mathcal{B})op(s)$ and
$v_{1}=op(\mathcal{B})(op(\hat{\chi}-1)+ \hat{R}_{m})v $.\\
Using the fact that supp$(1-\hat{\chi})\cap \mbox{supp}(\chi^{0})=
\varnothing$, $ \hat{R}_{m}\in \mathcal{T}\mathcal{O}^{-m-1}$, for $
m$ large and estimate (\ref{f1}), we show
\begin{equation}\label{m1}
\|F_{1}\|^{2}\leq C\left( \norm{P(x,D,\mu)
u}_{L^{2}}^{2}+\norm{u}_{1,\mu}^{2}\right).
\end{equation}
Using the fact that supp$(1-\hat{\chi})\cap \mbox{supp}(\chi^{0})=
\varnothing$, $\hat{R}_{m}\in \mathcal{T}\mathcal{O}^{-m-1}$, for $
m$ large and the trace formula (\ref{n3}), we show
\begin{equation}\label{m2}
\mu\abs{v_{1}}^{2}\leq C\left(\mu^{-2}
\abs{u}_{1,0,\mu}^{2}+\norm{u}_{1,\mu}^{2}\right).
\end{equation}
 Here we need to recall an argument shown in
Taylor \cite {Ta} given by this lemma
\begin{lemm}\label{l3}
Let $v$ solves the system $$ \frac{\partial}{\partial y}v=Gv+Av$$
where $ G= \left( \begin{array}{cc} E&\\
&F\end{array}\right)$ and $A$ are pseudo-differential operators of
order $1$ and $0$, respectively. We suppose that the symbols of $E$
and $F$ are two square matrices and have disjoint sets of
eigenvalues. Then there exists a pseudo-differential operator $K$ of
order $-1$ such that $w=(I+K)v$ satisfies
$$ \frac{\partial}{\partial y}w= Gw + \left( \begin{array}{cccc}
\alpha_{1}&\\
&\alpha_{2}\end{array}
 \right)w+ R_{1}w+R_{2}v$$
where $\alpha_{j}$ and $R_{j}$, $j=1,2$ are pseudo-differential
operators of order $0$ and $- \infty$, respectively.
\end{lemm}
\par By this argument, there exists a pseudo-differential operator
$K(x, D_{x'}, \mu)$ of order $-1$ such that the boundary problem
(\ref{f4}) is reduced to the following
\begin{equation}
\left\{\begin{array}{ll}\label{z5}
D_{x_{n}}w-op( \mathcal{H})w=\tilde{F}&\mbox{in}\,\{x_{n}>0\},\\\\
op(\tilde{\mathcal{B}})w=
(\frac{1}{\mu}\Lambda\tilde{e}_{1},\tilde{e}_{2}) +v_{1}+
v_{2}&\mbox{on}\,\{x_{n}=0\},
\end{array}
\right.
\end{equation}
where $w=(I+K)V$, $\tilde{F}=(I+K)F_{1}$, $op(\mathcal{H})$ is a
tangential of order $1$ with principal symbol
$\mathcal{H}=\mbox{diag}(\mathcal{H}^{-}, \mathcal{H}^{+})$ and
$-\mbox{Im}(\mathcal{H}^{-})\geq C\Lambda$,
$op(\tilde{\mathcal{B}})=op(\mathcal{B}_{1})(I+K')$ with $ K'$ is
such that $(I+K')(I+K)= Id+R'_{m}$ ($ R'_{m}\in \mathcal{O}^{-m-1}$,
for $m$ large) and $ v_{2}= op(\mathcal{B}_{1})R'_{m}V$.\\
According to (\ref{m1}), we have
\begin{equation}\label{f5}
\|\tilde{F}\|^{2}\leq C\left( \norm{P(x,D,\mu)
u}_{L^{2}}^{2}+\norm{u}_{1,\mu}^{2}\right).
\end{equation}
Using the fact that $ R'_{m}\in \mathcal{O}^{-m-1}$, for $ m$ large,
the trace formula (\ref{n3}) and estimates (\ref{f2}), (\ref{f3})
and (\ref{m2}), we show
\begin{equation}\label{f6}
\mu\abs{op(\tilde{\mathcal{B}})w}^{2}\leq C\left( \frac{1}{\mu}
\abs{e_{1}}_{1}^{2}+\mu \abs{e_{2}}^{2}+\mu^{-2}
\abs{u}_{1,0,\mu}^{2}+\norm{u}_{1,\mu}^{2}\right).
\end{equation}
\begin{lemm}\label{l4}
Let $\mathcal{R}=\mbox{diag}(-\rho \mbox{Id}_{2},0)$, $\rho>0$. Then
there exists $ C>0$ such that
\begin{enumerate}
  \item  $\mbox{Im}(\mathcal{RH})=\mbox{diag}\left(e(x,\xi',\mu),
  0\right)$, with $e(x,\xi',\mu)=-\rho Im(\mathcal{H}^{-})$,
  \item  $e(x,\xi',\mu)\geq C\Lambda$ in \mbox{supp} $(\chi^{0})$,
  \item  $-\mathcal{R}+ \tilde{\mathcal{B}}^{\star}\tilde{\mathcal{B}}\geq
  C.\mbox{Id}$ on $\{x_{n}=0\}\cap \mbox{supp}\,(\chi^{0})$.
\end{enumerate}
\end{lemm}
{\bf Proof}\\ Denote the principal symbol $\tilde{\mathcal{B}}$ of
the boundary operator $op (\tilde{\mathcal{B}})$ by $
\left(\tilde{\mathcal{B}}^{-},\tilde{\mathcal{B}}^{+} \right)$ where
$ \tilde{\mathcal{B}}^{+}$ is the restriction of $
\tilde{\mathcal{B}}$ to subspace generated by $\left(
s_{1}^{+},s_{2}^{+}\right)$. We begin by proving that
$\tilde{\mathcal{B}}^{+}$ is an isomorphism. Denote $$ w_{1}=^t(1,0)
\quad\mbox{and}\quad w_{2}=^t(0,1).$$ Then
$$ \left\{\begin{array}{l}
s_{1}^{+}=\left(w_{1}, z_{1}^{+}\Lambda^{-1}w_{1}\right)\\\\
s_{2}^{+}=\left(w_{2}, z_{2}^{+}\Lambda^{-1}w_{2}\right)
\end{array} \right.$$
are eigenvectors of $z_{1}^{+}$ and $z_{2}^{+}$. We have
$\tilde{\mathcal{B}}^{+}=(\mathcal{B}_{0}+\frac{1}{\mu}r_{0})(s_{1}^{+}\,s_{2}^{+})= \mathcal{B}_{0}^{+}+ \frac{1}{\mu}r_{0}^{+} $. To proof
that $\tilde{\mathcal{B}}^{+}$ is an isomorphism it suffices, for
large $\mu$, to proof that $\mathcal{B}_{0}^{+} $ is an isomorphism.
Following (\ref{e9}), we obtain
$$ \mathcal{B}_{0}^{+}=\left(
                         \begin{array}{cc}
                           0 & -i \\
                           \Lambda^{-1}i\alpha_{1} & \Lambda^{-1}i\alpha_{2}\\
                         \end{array}
                       \right).
$$
Then
$$ \mbox{det}(\mathcal{B}_{0}^{+})=- \Lambda^{-1}\alpha_{1}. $$
If we suppose that $ \mbox{det}(\mathcal{B}_{0}^{+})=0$, we obtain $
\alpha_{1}=0$ and then $\alpha_{1}^{2}=0$. Following (\ref{e8}), we
obtain
$$
q_{1}=0 \quad\mbox{and}\quad \left(\mu \frac{\partial
\varphi_{1}}{\partial x_{n}}\right)^{2}+q_{2,1}=0.
$$ Combining with the fact that
$q_{2,1}+\frac{q_{1}^{2}}{\left( \partial \varphi_{1}/\partial x_{n}
\right)^{2} }=0$, we obtain $\left( \mu \frac{\partial
\varphi_{1}}{\partial x_{n}}\right)^{2}=0$, that is impossible
because following (\ref{e10}), we have $ \left(\frac{\partial
\varphi_{1}}{\partial x_{n}}\right)^{2}\neq 0$ and following (
\ref{e5}), we have $ \mu \neq 0$. We deduce that
$\tilde{\mathcal{B}}^{+}$ is an isomorphism.\\ Let us show the Lemma
\ref{l4}. We have
\begin{equation}\label{z6}
\mbox{Im}(\mathcal{RH})=\mbox{diag}\left(-\rho\,\mbox{Im}(\mathcal{H}^{-}),0\right)=
\mbox{diag}\left(e(x,\xi',\mu),0\right),
\end{equation}
where $ e(x,\xi',\mu)=-\rho\,\mbox{Im}(\mathcal{H}^{-})\geq C
\Lambda$, $C>0$. It remains to proof 3. \\ Let $w=(w^{-},w^{+})\in
\C^{4}= \C^{2}\oplus\C^{2}$. Then we have $ \tilde{\mathcal{B}}w=
\tilde{\mathcal{B}}^{-}w^{-}+\tilde{\mathcal{B}}^{+}w^{+}$. Since
$\tilde{\mathcal{B}}^{+} $ is an isomorphism, then there exists a
constant $ C>0$ such that
$$ \abs{\tilde{\mathcal{B}}^{+}w^{+}}^{2}\geq C \abs{w^{+}}^{2}.$$
Therefore, we have $$\abs{w^{+}}^{2}\leq C\left(
\abs{\tilde{\mathcal{B}}w }^{2}+\abs{w^{-}}^{2}\right).$$ We deduce
$$ -(\mathcal{R}w,w)=\rho \abs{w^{-}}^{2}\geq \frac{1}{C}\abs{w^{+}}^{2}+(\rho-1)
\abs{w^{-}}^{2}-\abs{\tilde{\mathcal{B}}w }^{2}. $$ Then, we obtain
the result, if $ \rho$ is large enough.\\
$~~~~~~~~~~~~~~~~~~~~~~~~~~~~~~~~~~~~~~~~~~~~~~~~~~~~~~~~
~~~~~~~~~~~~~~~~~~~~~~~~~~~~~~~~~~~~~~~~~~~~~~~~~~~~\square$
\subsubsection{Proof of proposition \ref{p2}}\hspace{5mm} We start by showing
(\ref{z1}). We have
$$\begin{array}{lll}
\norm{P_{1}(x,D,\mu)u_{0}}^{2}&=&\norm{(\mbox{Re}P_{1})u_{0}}^{2}+\norm{(\mbox{Im}P_{1})u_{0}}^{2}\\&+&
i\Bigg[ \Big((\mbox{Im}P_{1})u_{0},
(\mbox{Re}P_{1})u_{0}\Big)-\Big((\mbox{Re}P_{1})u_{0},
(\mbox{Im}P_{1})u_{0}\Big) \Bigg].\end{array}$$ By integration by
parts we find
$$\norm{P_{1}(x,D,\mu)u_{0}}^{2}=\norm{(\mbox{Re}P_{1})u_{0}}^{2}+\norm{(\mbox{Im}P_{1})u_{0}}^{2}
+i\Big( \left[
\mbox{Re}P_{1},\mbox{Im}P_{1}\right]u_{0},u_{0}\Big)+\mu
Q_{0}(u_{0}),$$ where
$$\left\{\begin{array}{lcl}Q_{0}(u_{0})&=&(-2\frac{\partial\varphi_{1}}{\partial x_{n}}D_{x_{n}}u_{0},D_{x_{n}}u_{0})_{0}+(op(r_{1})u_{0},
D_{x_{n}}u_{0})_{0}\\\\&+&(op(r'_{1})D_{x_{n}}u_{0},u_{0})_{0}+(op(r_{2})u_{0},u_{0})_{0}+\mu
(\frac{\partial\varphi_{1}}{\partial x_{n}}
u_{0},u_{0} )_{0},\\\\
&&r_{1}=r'_{1}=2q_{1,1},\quad\quad
r_{2}=-2\frac{\partial\varphi_{1}}{\partial x_{n}}q_{2,1}.
\end{array}
\right.$$ Then we have
$$\abs{Q_{0}(u_{0})}^{2}\leq C\abs{u_{0}}_{1,0,\mu}^{2}.$$
We show the same thing for $P_{2}(x,D,\mu)v_{0}$. In addition we
know that the principal symbol of the operator $ [ \mbox{Re}P_{j},
\mbox{Im}P_{j}]$, $j=1,2$, is given by $
\{\mbox{Re}P_{j},\mbox{Im}P_{j}\}$. Proceeding like Lebeau and
Robbiano in paragraph 3 in \cite{LR}, we obtain (\ref{z1}).\par It
remains to prove (\ref{z2}). Following Lemma \ref{l4}, let
$G(x_{n})=d/dx_{n}(op(\mathcal{R})w,w)_{L^{2}(\R^{n-1})}$.\\ Using $
D_{x_{n}}w-op(\mathcal{H})=\tilde{F}$, we obtain
$$G(x_{n})=-2\,\mbox{Im}(op(\mathcal{R})\tilde{F},w)-2\,\mbox{Im}(op(\mathcal{R})op(\mathcal{H})w,w).
$$The integration in the normal direction gives
\begin{equation}\label{z7}
   (op(\mathcal{R})w,w)_{0}=
   \dint_{0}^{\infty}\mbox{Im}(op(\mathcal{R})op(\mathcal{H})w,w)dx_{n}+
    2\dint_{0}^{\infty}\mbox{Im}(op(\mathcal{R})\tilde{F},w)dx_{n}.
\end{equation}
From Lemma \ref{l4} and the G{\aa}rding inequality, we obtain, for
$\mu$ large,
\begin{equation}\label{z8}
\mbox{Im}(op(\mathcal{R})op(\mathcal{H})w,w)\geq
C\abs{w^{-}}_{\frac{1}{2}}^{2},
\end{equation}
moreover we have for all $\epsilon>0$
\begin{equation}\label{z9}
\int_{0}^{\infty}\abs{(op(\mathcal{R})\tilde{F},w)}dx_{n}\leq
\epsilon
C\mu\norm{w^{-}}^{2}+\frac{C_{\epsilon}}{\mu}\|\tilde{F}\|^{2}.
\end{equation}
Applying Lemma \ref{l4} and the G{\aa}rding inequality, we obtain,
for $\mu$ large,
\begin{equation}\label{z10}
    -(op(\mathcal{R})w,w)+|op(\tilde{\mathcal{B}})w|^{2}\geq C\abs{w}^{2}.
\end{equation}
 Combining (\ref{z10}), (\ref{z9}), (\ref{z8}) and (\ref{z7}),
we get
\begin{equation}\label{z11}
    C\abs{w^{-}}_{\frac{1}{2}}^{2}+C\abs{w}^{2}\leq
    \frac{C}{\mu}\|\tilde{F}\|^{2}+|op(\tilde{\mathcal{B}})w|^{2}.
\end{equation}
Then
$$ \mu \abs{w}^{2}\leq C\|\tilde{F}\|^{2}+ \mu |op(\tilde{\mathcal{B}})w|^{2}.$$
Recalling that $w=(I+K)V$, $V=op(n)v$, $v=^t\left(\langle
D',\mu\rangle \tilde{u},D_{x_{n}}\tilde{u}\right)$ and $
\tilde{u}=op(\chi^{0})u$ and using estimates (\ref{f5}) and (\ref{f6}), we prove (\ref{z2}).\\
$~~~~~~~~~~~~~~~~~~~~~~~~~~~~~~~~~~~~~~~~~~~~~~~~~~~~~~~~~~~~~~~~~~~~~~~~~~~~~~~~~~~~~~~~~~~~~~~~~~~~~~~~~~~~\square$
\subsection{Estimate in $\mathcal{E}_{1}^{-}$} This part is devoted
to estimate in region $ \mathcal{E}_{1}^{-}$.\\ Let
$\chi^{-}(x,\xi',\mu)\in \mathcal{T}S_{\mu}^{0}$ equal to $1$ in $
\mathcal{E}_{1}^{-}$ and such that in the support of $ \chi^{-}$ we
have $q_{2,1}+\frac{q_{1}^{2}}{(\partial \varphi_{1}/\partial
x_{n})^{2}} \leq -\delta < 0$. Then we have the following partial
estimate.
\begin{prop}\label{p3}
There exists constants $C >0$ and $\mu_{0}>0$ such that for any
$\mu\geq\mu_{0}$ we have the following estimate
\begin{equation}\label{k1}
\mu\norm{op(\chi^{-})u}_{1,\mu}^{2}\leq
C\left(\norm{P(x,D,\mu)u}^{2}+
 \mu\abs{u}_{1,0,\mu}^{2}+\norm{u}_{1,\mu}^{2} \right),
\end{equation}
for any $u\in C_{0}^{\infty}(\overline{\Omega}_{2})$.\par If we assume
moreover that $\frac{\partial \varphi_{1}}{\partial x_{n}}>0$ then
we have
\begin{equation}\label{k2}
\mu\abs{op(\chi^{-})u_{0}}_{1,0,\mu}^{2}\leq
C\left(\norm{P(x,D,\mu)u}^{2}+
 \mu^{-2}\abs{u}_{1,0,\mu}^{2}+\norm{u}_{1,\mu}^{2} \right)
\end{equation}
for any $u=(u_{0},v_{0})\in C_{0}^{\infty}(\overline{\Omega}_{2})$.
\end{prop}
{\bf Proof.}\\
Let $\tilde{u}=op(\chi^{-})u= (op(\chi^{-})u_{0},op(\chi^{-})v_{0}
)=(\tilde{u}_{0},\tilde{v}_{0})$.\\
In this region we have not a priori information for the roots of
$p_{2}(x,\xi,\mu)$. Using the same technique of the proof of
(\ref{z1}), we obtain
\begin{equation}\label{k3}
\mu\norm{op(\chi^{-})v_{0}}_{1,\mu}^{2}\leq
C\left(\norm{P(x,D,\mu)v_{0}}^{2}+
 \mu\abs{v_{0}}_{1,0,\mu}^{2}+\norm{v_{0}}_{1,\mu}^{2} \right)
\end{equation}
In supp$(\chi^{-})$ the two roots $ z_{1}^{\pm}$ of
$p_{1}(x,\xi,\mu)$ are in the half-plane $ Im\xi_{n}<0$. Then we can
use the Calderon projectors. By the same way that the proof of
(\ref{e17}) and using the fact that the operators $t_{0,1}$ and $
t_{1,1}$ vanish in $ x_{n}>0$ (because the roots are in
$Im\xi_{n}<0$, see (\ref{n5})), the counterpart of (\ref{n2}) is
then
\begin{equation}\label{g1}
\tilde{u}_{0}=E\underline{\tilde{f}}_{1}+w_{1,1}+w_{2,1},\quad\mbox{for}\,
x_{n}>0.
\end{equation}
We then obtain (see proof of (\ref{e17}))
\begin{equation}\label{g2}
\mu^{2}\norm{op(\chi^{-})u_{0}}_{1,\mu}^{2}\leq C\left(
\norm{P_{1}(x,D,\mu)u_{0}}^{2}+
 \mu\abs{u_{0}}_{1,0,\mu}^{2}+\norm{u_{0}}_{1,\mu}^{2}\right).
\end{equation}
Combining (\ref{k3}) and (\ref{g2}), we obtain (\ref{k1}).\\ It
remains to proof (\ref{k2}). We take the trace at $x_{n}=0^{+}$ of
(\ref{g1}),
$$ \gamma_{0}(\tilde{u}_{0})=w_{0,1}=
\gamma_{0}( E\underline{\tilde{f}}_{1}+w_{1,1}+w_{2,1}),$$ which, by
the counterpart of (\ref{n9}), gives
\begin{equation}\label{g3}
\mu\abs{\gamma_{0}(\tilde{u}_{0})}_{1}^{2}\leq C\left(
\norm{P_{1}(x,D,\mu)u_{0}}^{2}+\norm{u_{0}}_{1,\mu}^{2}+\mu^{-2}\abs{u_{0}}_{1,0,\mu}^{2}\right).
\end{equation}
From (\ref{g1}) we also have
$$ D_{x_{n}}\tilde{u}_{0}=D_{x_{n}}
E\underline{\tilde{f}}_{1}+D_{x_{n}}w_{1,1}+D_{x_{n}}w_{2,1},
\quad\mbox{for}\, x_{n}>0.$$ We take the trace at $x_{n}=0^{+}$ and
obtain
$$ \gamma_{1}(\tilde{u}_{0})=\gamma_{0}(D_{x_{n}}(E\underline{\tilde{f}}_{1}+w_{1,1}+w_{2,1} )).$$
Using the trace formula (\ref{n3}), we obtain
$$ \abs{\gamma_{1}(\tilde{u}_{0})}^{2}\leq C\mu^{-1} \norm{D_{x_{n}}(E\underline{\tilde{f}}_{1}+w_{1,1}+w_{2,1} ) }^{2}_{1,\mu}
\leq C\mu^{-1}
\norm{E\underline{\tilde{f}}_{1}+w_{1,1}+w_{2,1}}_{2,\mu}^{2}$$ and,
by the counterpart of (\ref{n0}), (\ref{n01}) and (\ref{n4}), this
yields
\begin{equation}\label{g4}
\mu\abs{\gamma_{1}(\tilde{u}_{0})}^{2}\leq C\left(
\norm{P_{1}(x,D,\mu)u_{0}}^{2}+\norm{u_{0}}_{1,\mu}^{2}+\mu^{-2}\abs{u_{0}}_{1,0,\mu}^{2}\right).
\end{equation}
Combining (\ref{g3}) and (\ref{g4}), we obtain
$$ \mu\abs{op(\chi^{-})u_{0}}_{1,0,\mu}^{2}\leq C\left(
\norm{P_{1}(x,D,\mu)u_{0}}^{2}+\norm{u_{0}}_{1,\mu}^{2}+\mu^{-2}\abs{u_{0}}_{1,0,\mu}^{2}\right).$$
Then we have
(\ref{k2}).$~~~~~~~~~~~~~~~~~~~~~~~~~~~~~~~~~~~~~~~~~~~~~~~~~~~~~~~~~~~~~~~~~~~~~~~~~~~~~~~~\square$
\subsection{End of the proof} We choose a partition of unity
$\chi^{+}+\chi^{0}+\chi^{-}=1$ such that $\chi^{+}$, $\chi^{0}$ and
$\chi^{-}$ satisfy the properties listed in proposition \ref{p1},
\ref{p2} and \ref{p3} respectively. We have $$
\norm{u}_{1,\mu}^{2}\leq
\norm{op(\chi^{+})u}_{1,\mu}^{2}+\norm{op(\chi^{0})u}_{1,\mu}^{2}+\norm{op(\chi^{-})u}_{1,\mu}^{2}.$$
Combining this inequality and (\ref{e17}), (\ref{z1}) and
(\ref{k1}), we obtain, for large $\mu$, the first estimate
(\ref{k4}) of Theorem \ref{t2}. i.e.
$$
\mu\norm{u}_{1,\mu}^{2}\leq C\left( \norm{P(x,D,\mu)u}^{2}+\mu
\abs{u}_{1,0,\mu}^{2}\right). $$
 It remains to estimate $\mu
\abs{u}_{1,0,\mu}^{2}$. We begin by giving an estimate of $\mu
\abs{u_{0}}_{1,0,\mu}^{2}$.\\ We have
$$ \abs{u_{0}}_{1,0,\mu}^{2}\leq \abs{op(\chi^{+})u_{0}}_{1,0,\mu}^{2}+\abs{op(\chi^{0})u_{0}}_{1,0,\mu}^{2}+
\abs{op(\chi^{-})u_{0}}_{1,0,\mu}^{2}, $$
$$\abs{op(\chi^{+})u_{0}}_{1,0,\mu}^{2}\leq \abs{op(\chi^{+})u}_{1,0,\mu}^{2}
$$ and $$ \abs{op(\chi^{0})u_{0}}_{1,0,\mu}^{2}\leq \abs{op(\chi^{0})u}_{1,0,\mu}^{2}.$$
Combining these inequalities, (\ref{e18}), (\ref{z2}), (\ref{k2})
and the fact that\\ $
\mu^{-2}\abs{u}_{1,0,\mu}^{2}=\mu^{-2}\abs{u_{0}}_{1,0,\mu}^{2}+\mu^{-2}\abs{v_{0}}_{1,0,\mu}^{2}
$, we obtain, for large $\mu$
\begin{equation}\label{k5}
\mu\abs{u_{0}}_{1,0,\mu}^{2}\leq
C\left(\norm{P(x,D,\mu)u}^{2}+\mu^{-1}\abs{op(b_{1})u}_{1}^{2}+\mu\abs{op(b_{2})u}^{2}+\mu^{-2}\abs{v_{0}}_{1,0,\mu}^{2}
+\norm{u}^{2}_{1,\mu}\right).
\end{equation}
For estimate $\mu \abs{v_{0}}_{1,0,\mu}^{2}$, we need to use the
transmission conditions given by (\ref{e16}). We have
$$op(b_{1})u=u_{0}|_{x_{n}=0}-i\mu v_{0}|_{x_{n}=0}\quad\mbox{on}\,
\left\{x_{n}=0\right\}.$$ Then
$$ \mu\abs{v_{0}}_{1}^{2}\leq C\left( \mu^{-1} \abs{u_{0}}_{1}^{2}+\mu^{-1}\abs{op(b_{1})u}_{1}^{2}\right).$$
Since we have $ \mu^{-1}\abs{u_{0}}_{1}^{2}\leq
\mu\abs{u_{0}}_{1,0,\mu}^{2}$. Then using (\ref{k5}), we obtain
\begin{equation}\label{k6}
\mu\abs{v_{0}}_{1}^{2}\leq
C\left(\norm{P(x,D,\mu)u}^{2}+\mu^{-1}\abs{op(b_{1})u}_{1}^{2}+\mu
\abs{op(b_{2})u}^{2}+\mu^{-2}\abs{v_{0}}_{1,0,\mu}^{2}+\norm{u}^{2}_{1,\mu}
\right).
\end{equation}
We have also
$$ op(b_{2})u=\left(D_{x_{n}}+i\mu \frac{\partial
\varphi_{1}}{\partial
x_{n}}\right)u_{0}|_{x_{n}=0}+\left(D_{x_{n}}+i\mu \frac{\partial
\varphi_{2}}{\partial
x_{n}}\right)v_{0}|_{x_{n}=0}\quad\mbox{on}\,\left\{x_{n}=0\right\}.$$
Then$$
 \mu\abs{D_{x_{n}}v_{0}}^{2}\leq C\left( \mu\abs{op(b_{2})u}^{2}+
 \mu \abs{D_{x_{n}}u_{0}}^{2}+\mu^{3}\abs{u_{0}}^{2}+\mu^{3}\abs{v_{0}}^{2}\right).$$
 Using the fact that $ \abs{u}_{k-1}\leq
 \mu^{-1}\abs{u}_{k}$, we obtain
$$\mu\abs{D_{x_{n}}v_{0}}^{2}\leq C\left( \mu\abs{op(b_{2})u}^{2}+
 \mu \abs{D_{x_{n}}u_{0}}^{2}+\mu\abs{u_{0}}_{1}^{2}+\mu\abs{v_{0}}_{1}^{2}\right).$$
Since we have
$\mu\abs{u_{0}}_{1,0,\mu}^{2}=\mu\abs{D_{x_{n}}u_{0}}^{2}+\mu\abs{u_{0}}_{1}^{2}$.
Then using (\ref{k5}) and (\ref{k6}), we obtain
\begin{equation}\label{k7}
\mu\abs{D_{x_{n}}v_{0}}^{2}\leq C\left(
\norm{P(x,D,\mu)u}^{2}+\mu^{-1}\abs{op(b_{1})u}_{1}^{2}+\mu
\abs{op(b_{2})u}^{2}+\mu^{-2}\abs{v_{0}}_{1,0,\mu}^{2}+
\norm{u}^{2}_{1,\mu}\right).
\end{equation}
Combining (\ref{k6}) and (\ref{k7}), we have
\begin{equation}\label{k8}
\mu\abs{v_{0}}_{1,0,\mu}^{2}\leq C\left(
\norm{P(x,D,\mu)u}^{2}+\mu^{-1}\abs{op(b_{1})u}_{1}^{2}+\mu
\abs{op(b_{2})u}^{2}+\norm{u}^{2}_{1,\mu}\right).
\end{equation}
Combining (\ref{k5}) and (\ref{k8}), we obtain
\begin{equation}\label{k9}
\mu\abs{u}_{1,0,\mu}^{2}\leq C\left(
\norm{P(x,D,\mu)u}^{2}+\mu^{-1}\abs{op(b_{1})u}_{1}^{2}+\mu
\abs{op(b_{2})u}^{2}+\norm{u}^{2}_{1,\mu}\right).
\end{equation}
Inserting (\ref{k9}) in (\ref{k4}) and for large $\mu$, we obtain
(\ref{tt2}).\\$~~~~~~~~~~~~~~~~~~~~~~~~~~~~~~~~~~~~~~~~~~~~~~~~~~~~~~~~~~~~~~~~~~~~~~~~~~~~~~~~~~~~~~~~~~~~~~~~~~~~~~~~~~\square$\\\\
\newpage
\section*{Appendix A}
This appendix is devoted to prove Lemma \ref{le1}. For this, we need
to distinguish two cases.
\begin{enumerate}
  \item {\bf {Inside $\mathcal{O}$}}\\
  To simplify the writing, we note $
  \norm{u}_{L^{2}(\mathcal{O})}=\norm{u}$.\\
  Let $\chi\in C_{0}^{\infty}(\mathcal{O})$. We have by integration
  by part
  $$ ((\triangle-i\mu)u, \chi^{2}u)=(-\nabla u,\chi^{2}\nabla u)-(\nabla u, \nabla (\chi^{2})u)
  -i\mu \norm{\chi u}^{2}.$$
    Then
  $$ \mu  \norm{\chi u}^{2}\leq C \left( \norm{f}\norm{\chi^{2}u}
  +\norm{\nabla u}^{2}+ \norm{\nabla u}\norm{\chi u}\right).$$
  Then
  $$\mu  \norm{\chi u}^{2}\leq C \left(\frac{1}{\epsilon}\norm{f}^{2}+\epsilon \norm{\chi^{2}u}
  +\norm{\nabla u}^{2}+ \frac{1}{\epsilon}\norm{\nabla u}^{2}+ \epsilon\norm{\chi u}^{2}\right) .$$
  Recalling that $\mu\geq1$, we have for $\epsilon$ small enough
 \begin{equation}\label{ap1}
  \norm{\chi u}^{2}\leq C \left( \norm{\nabla u}^{2} +\norm{f}^{2} \right).
 \end{equation}
  Hence the result inside $\mathcal{O}$.
 \item {\bf{ In the neighborhood of the boundary}}\\
 Let $x=(x',x_{n})\in \R^{n-1} \times \R$. Then $$\partial \mathcal{O}=\{x\in \R^{n},\,\,
 x_{n}=0\}.$$
Let $ \epsilon>0$ such that $0<x_{n}<\epsilon$. Then we have$$
 u(x',\epsilon)-u(x',x_{n})=\int_{x_{n}}^{\epsilon}\partial_{x_{n}}u(x',\sigma)d\sigma .$$
 Then $$ \abs{u(x',x_{n}) }^{2}\leq 2 \abs{u(x',\epsilon)}^{2}+ 2\left(\int_{x_{n}}^{\epsilon}\abs{\partial_{x_{n}}u(x',\sigma)}d\sigma \right)^{2}.$$
 Using the Cauchy Schwartz inequality, we obtain
$$
\abs{u(x',x_{n}) }^{2}\leq 2 \abs{u(x',\epsilon)}^{2}+2
\epsilon^{2}\int_{0}^{\epsilon}
\abs{\partial_{x_{n}}u(x',x_{n})}^{2}d x_{n}.
$$
Integrating with regard to $x'$, we obtain
\begin{equation}\label{ap2}
\int_{\abs{x'}<\epsilon}\abs{u(x',x_{n}) }^{2}dx'\leq 2
\int_{\abs{x'}<\epsilon}\abs{u(x',\epsilon)}^{2}dx'+2 \epsilon^{2}
\int_{\abs{x'}<\epsilon,\,\abs{x_{n}}<\epsilon}\left(
\abs{\partial_{x_{n}}u(x',x_{n})}^{2}d x_{n}\right)dx'.
\end{equation}
Using the trace Theorem, we have
\begin{equation}\label{ap3}
\int_{\abs{x'}<\epsilon}\abs{u(x',\epsilon)}^{2}dx'\leq C
\int_{\abs{x'}<2\epsilon,\, \abs{x_{n}-\epsilon}<\frac{\epsilon}{2}}
(\abs{u(x)}^{2}+\abs{\nabla u(x)}^{2})dx.
\end{equation}
Now we need to introduce the following cut-off functions
$$ \chi_{1}(x)=\left \{\begin{array}{lcl}
1&\mbox{if}& 0<x_{n}<\frac{\epsilon}{2},\\\\
0&\mbox{if}& x_{n>}\epsilon
\end{array}\right. $$
and
$$ \chi_{2}(x)=\left \{\begin{array}{lcl}
1&\mbox{if}& \frac{\epsilon}{2}<x_{n}<\frac{3 \epsilon}{2},\\\\
0&\mbox{if}& x_{n}<\frac{\epsilon}{4}, \,\,x_{n>}2\epsilon .
\end{array}\right. $$
Combining (\ref{ap2}) and (\ref{ap3}), we obtain for $\epsilon$
small enough
\begin{equation}\label{ap4}
\norm{\chi_{1}u}^{2}\leq C\left(\norm{\chi_{2}u}^{2}+ \norm{\nabla
u}^{2} \right).
\end{equation}
Since following (\ref{ap1}), we have
$$\norm{\chi_{2}u}^{2}\leq C\left( \norm{f}^{2}+ \norm{\nabla
u}^{2}\right). $$ Inserting in (\ref{ap4}), we obtain
\begin{equation}\label{ap5}
\norm{\chi_{1}u}^{2}\leq C\left( \norm{f}^{2}+ \norm{\nabla
u}^{2}\right).
\end{equation}
Hence the result in the neighborhood of the boundary.
\end{enumerate}
Following (\ref{ap1}), we can write
\begin{equation}\label{ap6}
\norm{(1-\chi_{1})u}^{2}\leq C\left( \norm{f}^{2}+ \norm{\nabla
u}^{2}\right).
\end{equation}
Adding (\ref{ap5}) and (\ref{ap6}), we obtain $$\norm{u}^{2}\leq
C\left( \norm{f}^{2}+ \norm{\nabla u}^{2}\right).$$ Hence the
result.
\newpage
\section*{Appendix B: Proof of Lemma \ref{le2}}
\hspace{5mm} This appendix is devoted to prove Lemma \ref{le2}.\\
Let $\chi\in C_{0}^{\infty}(\R^{n})$ such that $\chi=1$ in the
support of $u$. We want to show that $ op(\Lambda^{s})e^{\mu
\varphi}\chi op(\Lambda^{-s})$ is bounded in $L^{2}$. Recalling that
for all $u$ and $v \in \mathcal{S}(\R^{n})$, we have
$$ \mathcal{F}(uv)(\xi')=(\frac{1}{2\pi})^{n-1} \mathcal{F}(u)\ast \mathcal{F}(v)(\xi'),\quad\quad \forall \xi'\in\R^{n-1}.$$
Then
$$\begin{array}{lcl} \mathcal{F}(op(\Lambda^{s})e^{\mu \varphi}\chi
op(\Lambda^{-s})v)(\xi',\mu)&=& \langle\xi',\mu\rangle^{s }
\mathcal{F}(e^{\mu \varphi}\chi
op(\Lambda^{-s})v)(\xi',\mu)\\\\
&=&(\dfrac{1}{2\pi})^{n-1}\langle\xi',\mu\rangle^{s}(g(\xi',\mu)\ast
\langle\xi',\mu\rangle^{-s}\mathcal{F}(v))(\xi',\mu),
\end{array}
$$where $ g(\xi',\mu)= \mathcal{F}(e^{\mu\varphi}\chi)(\xi',\mu)$.
Then we have
$$ \mathcal{F}(op(\Lambda^{s})e^{\mu \varphi}\chi
op(\Lambda^{-s})v)(\xi',\mu)=\int g(\xi'-\eta',\mu)\langle
\xi',\mu\rangle^{s}\langle \eta',\mu\rangle^{-s}
\mathcal{F}(v)(\eta',\mu)d\eta'.$$ Let
$k(\xi',\eta')=g(\xi'-\eta',\mu) \langle \xi',\mu\rangle^{s}\langle
\eta',\mu\rangle^{-s}$. Our goal is to show that $\int K(\xi',\eta')
\mathcal{F}(v)(\eta',\mu)d\eta'$ is bounded in $L^{2}$. To do it, we
will use Lemma of Schur. It suffices to prove that there exists
$M>0$ and $N>0$ such that
$$ \int \abs{K(\xi',\eta')}d\xi'\leq M\quad\quad
\mbox{and}\quad\quad  \int \abs{K(\xi',\eta')}d\eta'\leq N .$$ In
the sequel, we suppose $s\geq 0$ (the case where $s<0$ is treated in
the same way).\\For $R>0$, we have
$$\begin{array}{lcl}
\langle \xi',\mu\rangle^{2R}g(\xi',\mu)&=&\dint \langle
\xi',\mu\rangle^{2R}e^{-i x'\xi'}\xi(x)e^{\mu\varphi(x)}dx'\\\\
&=&\dint(1-\Delta+\mu^{2})^{R}(e^{-i x'\xi'})\chi(x)e^{\mu\varphi(x)}dx'\\\\
&=&\dint e^{-i
x'\xi'}(1-\Delta+\mu^{2})^{R}(\chi(x)e^{\mu\varphi(x)})dx'.
\end{array}$$
Then there exists $C>0$, such that
\begin{equation}\label{ab}
\abs{\langle \xi',\mu\rangle^{2R}g(\xi',\mu)}\leq Ce^{C\mu}.
\end{equation}
Moreover, we can write
$$ \int\abs{K(\xi',\eta')}d\xi'= \int\left|{ g(\xi'-\eta',\mu) \langle\xi'-\eta',\mu \rangle^{2R}
\frac{\langle \xi',\mu\rangle^{s}\langle
\eta',\mu\rangle^{-s}}{\langle\xi'-\eta',\mu
\rangle^{2R}}}\right|d\xi'.$$ Using (\ref{ab}), we obtain
$$\int\abs{K(\xi',\eta')}d\xi'\leq C e^{C\mu}\int \frac{\langle \xi',\mu\rangle^{s}\langle
\eta',\mu\rangle^{-s}}{\langle\xi'-\eta',\mu \rangle^{2R}}d\xi'.$$
Since $$\int \frac{\langle \xi',\mu\rangle^{s}\langle
\eta',\mu\rangle^{-s}}{\langle\xi'-\eta',\mu \rangle^{2R}}d\xi'=
\int_{\abs{\xi'}\leq \frac{1}{\epsilon}\abs{\eta'}} \frac{\langle
\xi',\mu\rangle^{s}\langle
\eta',\mu\rangle^{-s}}{\langle\xi'-\eta',\mu \rangle^{2R}}d\xi'+
\int_{\abs{\eta'}\leq \epsilon\abs{\xi'}} \frac{\langle
\xi',\mu\rangle^{s}\langle+
\eta',\mu\rangle^{-s}}{\langle\xi'-\eta',\mu
\rangle^{2R}}d\xi',\quad \epsilon >0. $$
 If $\abs{\xi'}\leq \frac{1}{\epsilon}\abs{\eta'}$,
 we have $$ \frac{\langle \xi',\mu\rangle^{s}\langle
\eta',\mu\rangle^{-s}}{\langle\xi'-\eta',\mu \rangle^{2R}}\leq C
\frac{\langle \eta',\mu\rangle^{s}\langle
\eta',\mu\rangle^{-s}}{\langle\xi'-\eta',\mu \rangle^{2R}}\leq
\frac{C}{\langle\xi'-\eta',\mu \rangle^{2R}}\quad\in\,
L^{1}\quad\mbox{if} \quad 2R>n-1.
$$
 If $\abs{\eta'}\leq
\epsilon\abs{\xi'} $, i.e $\langle\xi'-\eta',\mu \rangle\geq \delta
\langle \xi',\mu\rangle$, $\delta>0$, we have
$$\frac{\langle \xi',\mu\rangle^{s}\langle
\eta',\mu\rangle^{-s}}{\langle\xi'-\eta',\mu \rangle^{2R}}\leq
\frac{C}{\langle\xi'-\eta',\mu \rangle^{2R-s}}\quad\in\,
L^{1}\quad\mbox{if} \quad 2R-s>n-1. $$ Then there exists $M>0$, such
that $$\int \abs{K(\xi',\eta')}d\xi'\leq M e^{C\mu}.$$ By the same
way, we show that there exists $N>0$, such that
$$\int \abs{K(\xi',\eta')}d\eta'\leq N e^{C\mu}.$$
Using Lemma of Schur, we have $ (op(\Lambda^{s})e^{\mu \varphi}\chi
op(\Lambda^{-s}))$ is bounded in $L^{2}$ and
$$ \norm{op(\Lambda^{s})e^{\mu \varphi}\chi
op(\Lambda^{-s})}_{ \mathcal{L}(L^{2})}\leq Ce^{C\mu}.$$ Applying in
$op(\Lambda^{s})u$, we obtain the result.
\bibliographystyle{plain}
\bibliography{biblio}
\end{document}